\documentclass[bj,preprint]{imsart}
\usepackage{latexsym,epsfig,amssymb,amsmath,amsthm,color,url,tikz,pdfsync}
\usetikzlibrary{arrows,positioning} 
\usetikzlibrary{shapes}

\RequirePackage[numbers,sort]{natbib}
\usepackage{graphicx,bbm}
\usepackage{inputenc}
\setcitestyle{numbers,multirow,import}
\usepackage{graphicx}
\allowdisplaybreaks
\numberwithin{equation}{section}
\numberwithin{figure}{section}

\startlocaldefs
\newtheorem{Theorem}{Theorem}[section]
\newtheorem{Lemma}[Theorem]{Lemma}
\newtheorem{Remark}[Theorem]{Remark}
\newtheorem{Example}[Theorem]{Example}
\newtheorem{Proposition}[Theorem]{Proposition}
\newtheorem{Definition}[Theorem]{Definition}
\newtheorem{Corollary}[Theorem]{Corollary}

\newcommand{\deltain}{\delta_{\text in}}
\newcommand{\deltaout}{\delta_{\text out}}

\newcommand{\bthe}{\begin{Theorem}}
\newcommand{\ethe}{\end{Theorem}}

\newcommand{\ble}{\begin{Lemma}}
\newcommand{\ele}{\end{Lemma}}

\newcommand{\bde}{\begin{Definition}}
\newcommand{\ede}{\end{Definition}}

\newcommand{\bco}{\begin{Corollary}}
\newcommand{\eco}{\end{Corollary}}

\newcommand{\bpr}{\begin{Proposition}}
\newcommand{\epr}{\end{Proposition}}

\newcommand{\brem}{\begin{Remark}}
\newcommand{\erem}{\end{Remark}}

\newcommand{\bexam}{\begin{Example}}
\newcommand{\eexam}{\end{Example}}

\newcommand{\beqq}{\begin{equation}}
\newcommand{\eeqq}{\end{equation}}

\newcommand{\beao}{\begin{eqnarray*}}
\newcommand{\eeao}{\end{eqnarray*}\noindent}

\newcommand{\beam}{\begin{eqnarray}}
\newcommand{\eeam}{\end{eqnarray}\noindent}

\newcommand{\barr}{\begin{array}}
\newcommand{\earr}{\end{array}}

\newcommand{\bproof}{\begin{proof}}
\newcommand{\eproof}{\end{proof}}

\newcommand{\sid}[1]{{\color{black} #1}}

\newcommand{\pin}{p^{\text{in}}}
\newcommand{\pout}{p^{\text{out}}}

\newcommand{\ain}{\iota_\text{in}}
\newcommand{\aout}{\iota_\text{out}}

\newcommand{\dd}{\mathrm{d}}
\newcommand{\PP}{\textbf{P}}
\newcommand{\EE}{\textbf{E}}
\newcommand{\ind}{\textbf{1}}

\newcommand{\convp}{\stackrel{P}{\longrightarrow}}

\newcommand{\convas}{\stackrel{\text{a.s.}}{\longrightarrow}}

\newcommand{\barV}{\overline{V}}
\newcommand{\barbV}{\overline{\mathbf{V}}}

\newcommand{\hatain}{\hat{\iota}_\text{in}}

\newcommand{\barbfI}{\overline{\mathbf{I}}}
\newcommand{\barbfO}{\overline{\mathbf{O}}}
\newcommand{\bfI}{\mathbf{I}}
\newcommand{\bfO}{\mathbf{O}}
\newcommand{\barI}{\overline{I}}
\newcommand{\barO}{\overline{O}}
\newcommand{\bfV}{\mathbf{V}}
\newcommand{\bfE}{\mathbf{E}}
\newcommand{\bfG}{\mathbf{G}}

\newcommand{\supy}{\sup_{y\ge 1}}

\newcommand{\argmin}{\operatornamewithlimits{argmin}}

\endlocaldefs

\begin{document}
\begin{frontmatter}
\title{A Directed Preferential Attachment Model with Poisson Measurement}
\runtitle{Poisson PA Model}

\begin{aug}
\author[A]{\fnms{Tiandong} \snm{Wang}\ead[label=e1]{twang@stat.tamu.edu}}
\and
\author[B]{\fnms{Sidney I.} \snm{Resnick}\ead[label=e2]{sir1@cornell.edu}}
\address[A]{Department of Statistics, Texas A\&M University, 
College Station, TX 77843, USA, \printead{e1}}

\address[B]{School of Operations Research and Information Engineering,
Cornell University, Ithaca, NY 14853, USA, \printead{e2}}
\end{aug}

\begin{abstract}
When modeling a directed social network, one choice is to use the traditional preferential attachment model, which generates power-law tail distributions.
In a traditional directed preferential attachment, every new edge is added sequentially into the network.
However, for real datasets, it is common to only have coarse timestamps available, which means several new edges are created at the same timestamp. 
Previous analyses on the evolution of social networks reveal that after reaching a stable phase, the growth of edge counts in a network follows a non-homogeneous Poisson process with a constant rate across the day but varying rates from day to day. Taking such empirical observations into account, we propose a modified preferential attachment model with Poisson measurement, 
and study its asymptotic behavior. This modified model is then fitted to real datasets, and 
we see it provides a better fit than the traditional one.
\end{abstract}


\end{frontmatter}




\section{Introduction}
Empirical evidence suggests the in- and out-degree distributions for nodes in
many social networks have Pareto-like tails (cf. \cite{kunegis:2013}).
A traditional preferential attachment (PA) model (\cite{bollobas:borgs:chayes:riordan:2003, krapivsky:redner:2001})
 theoretically generates a network that exhibits such heavy-tailed properties
under the intuitive assumption that nodes with large degrees tend to attract more edges than those with small degrees. 
For these reasons, the traditional PA model has attracted a great amount of attention in the modeling of social networks.

However, sometimes simple assumptions do not match with what we have observed from real datasets.
For example, in a traditional directed PA setup (cf. \cite{wan:wang:davis:resnick:2017, krapivsky:redner:2001, bollobas:borgs:chayes:riordan:2003, resnick:samorodnitsky:towsley:davis:willis:wan:2016}),
every new edge is added sequentially, annotated with a unique timestamp upon its creation.
But in a lot of real examples (e.g. the second dataset in Section~\ref{sec:intro_data}), timestamp information is coarse and 
it is possible to have more than one edge created at one single timestamp.

Here we first discuss two real data examples, \emph{Facebook wall posts} and \emph{Slashdot reply network}, from which
we summarize important features that are not captured by the traditional PA model. 
Based on these observed features, we propose a modified directed PA model in Section~\ref{subsec:POPA} and study its theoretical properties.

\subsection{Data examples}\label{sec:intro_data}
\subsubsection{Facebook wall posts}\label{subsec:intro_fb}
In \cite{wang:resnick:2019b}, several geographically concentrated networks have been studied and a common pattern 
in the growth of a network has been observed. Empirical findings suggest that the start-up phase of the growth of edge counts in a regional
network can be modeled by a self-exciting point process. After the start-up phase ends,
the growth of the edge counts can be modeled instead by a non-homogeneous Poisson process (NHPP) with a
constant rate across the day but varying rates from day to day, plus a nightly inactive
period when local users are expected to be asleep.

One particular example considered in \cite{wang:resnick:2019b} is the Facebook wall post data for users in New Orleans
(available at \url{http://konect.uni-koblenz.de/networks/facebook-wosn-wall}),
and timestamps generated during the expected daily sleeping hours 1-8 AM U.S. Central Time have been excluded from our modeling. 
After the start-up phase of this regional Facebook network ends, we model the 
edge creation process by an NHPP with constant rates within a day but varying rates from day to day.
We then model the node creation process by another NHPP.
Applying a Kolmogorov-Smirnov (KS) test to the node creation process shows the plausibility of fitting an NHPP, and
no significant evidence flagging the dependence among residuals has been detected using Ljung-Box tests.
So it is also reasonable to view the node creation process as thinning the edge creation process.

\begin{figure}
\centering
\includegraphics[scale=.5]{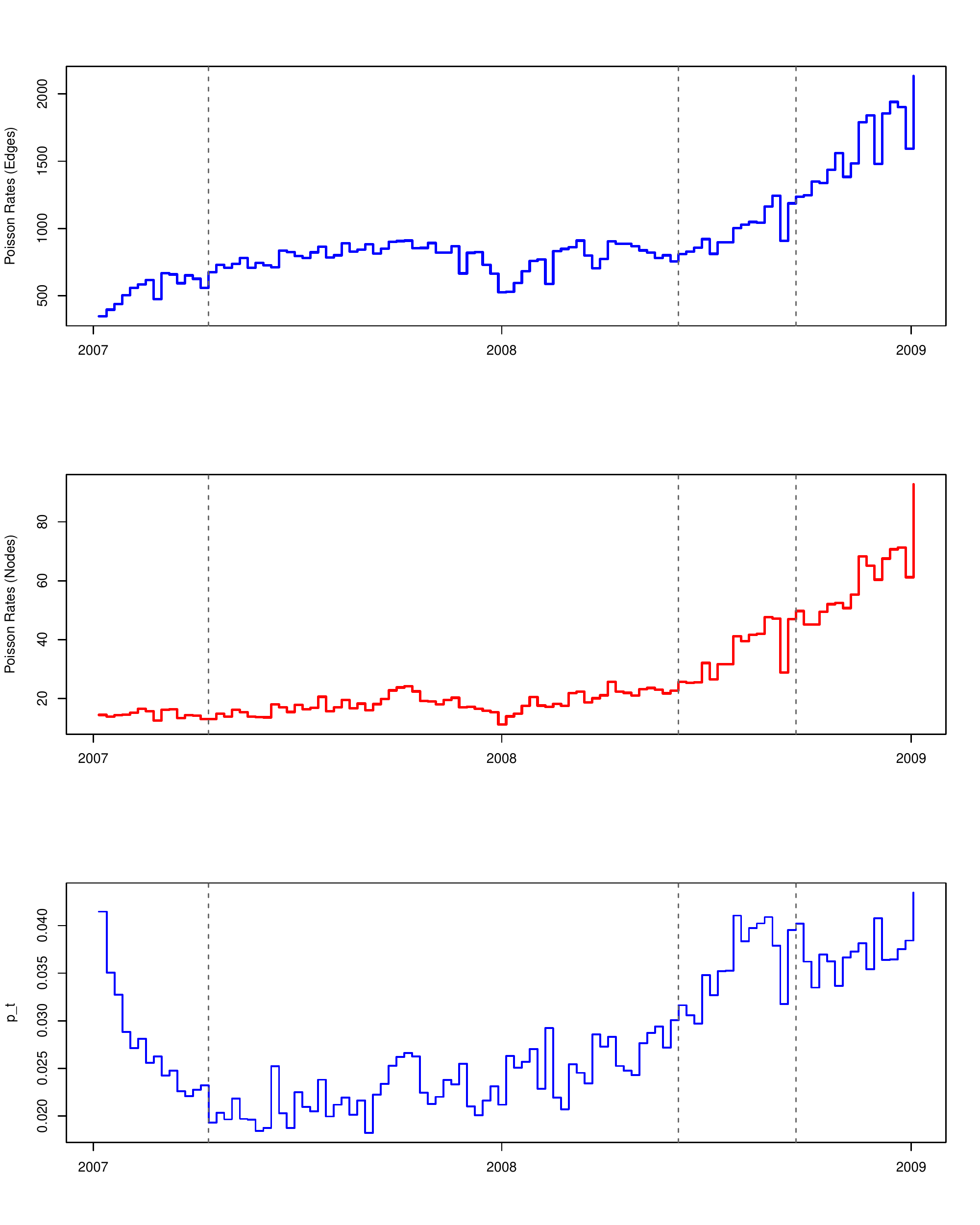}
\caption{Estimated daily Poisson rates of the edge (top), node (middle) creation process and the ratio of the two Poisson rates $\widehat p_t$ (bottom). 
Within the second time segment, all three quantities remain relatively stable.}\label{fig:pois}
\end{figure}

We calculate the daily Poisson rate estimates for both edge and node creation processes, 
and average them over non-overlapping weekly intervals.
A graphical illustration is given in the top and middle panels of Figure~\ref{fig:pois}.
The bottom panel of Figure~\ref{fig:pois} reports the ratio 
\[
\widehat p_t = \frac{\text{Daily Poisson Rate Estimates for the Node Creation Process}}{\text{Daily Poisson Rate Estimates for the Edge Creation Process}}.
\]
We use the \verb6breakpoints6 function in R's \verb6strucchange6 package to identify change points in
the daily Poisson rate estimates of the edge generation process, and in Figure~\ref{fig:pois} they are denoted by the grey vertical lines.
Within the second time segment, we see that all three quantities in Figure~\ref{fig:pois} remain relatively stable.

\subsubsection{Slashdot}\label{subsec:intro_sd}
Another data example is the reply network of the technology website, \emph{Slashdot}, which is available at
\url{http://konect.cc/networks/slashdot-threads/}. 
In this network, nodes correspond to different users, and directed edges represent the replies which start from the responding user. 

Although the dataset contains timestamp information,
the minimum time elapsed between two adjacent timestamps is counted in minutes. In other words,
the timestamp is coarse, and several edge creation events may happen at the same timestamp. 
This raises difficulties in model fitting since standard methods like MLE \citep{wan:wang:davis:resnick:2017} require knowing the exact
evolution history of the network, and assume each edge is created at a unique timestamp. 
Unlike Facebook, a node labeled as a smaller number in the Slashdot data
is not necessarily created at an earlier time in the network. 
So inferring the order of edges created in the network through node labels does not work, either.
The coarse timestamp may also lead to lags in updating the configuration of the network. When a new 
edge is created, the configuration is not updated until a later timestamp.
Such delay cannot be modeled by a traditional PA model where new edges are added sequentially.

Another phenomenon observed is that both the Facebook and Slashdot networks  display empirical in- and out-degree distributions with power-law tails 
(see Figure~\ref{fig:fb_PoisPA} and \ref{fig:sd_PoisPA} for example).

\subsection{Goals}
Motivated by the observations from the two datasets, we summarize that a modified network model is necessary and
it must:
\begin{enumerate}
\item[(i)] Allow the creation of a Poisson number of nodes and edges at each step of the network evolution.
\item[(ii)] Capture the possibility that the timestamp information may be only coarsely observed.
\item[(iii)] Generate in- and out-degree distributions with power-law tails.
\end{enumerate}
Note that the traditional PA model fails to capture the first two observations. In this paper, we modify the PA assumptions by taking the
first two findings into account, and study the asymptotic properties of the modified model such that
the third feature is guaranteed.

\medskip

The rest of the paper is organized as follows. In Section~\ref{sec:prelim}, we give the description of our modified PA network model and compare it with 
the traditional directed PA model. The formal constructions for the traditional PA and the modified PA model with Poisson measurement are
developed in Section \ref{sec:embed_PA} and \ref{sec:POPA}, respectively.
Relevant convergence results for the modified PA model are given in Section~\ref{subsec:POPA_deg}, and we include discussions on model fitting in Section~\ref{sec:fit}.
Additional comments are given in Section~\ref{sec:comments}, and all
proofs are collected in Section~\ref{sec:proof}.
In fact, based on the proof machinery in Section~\ref{sec:proof}, we can also
relax the Poisson assumption on the edge creation process to any iid non-negative random variables with finite first inverse moment.

\section{Description of the Two PA Models}\label{sec:prelim}
Taking the observations in Section \ref{sec:intro_data} into account, in this section, we describe a modified linear PA model by adding a Poisson number of 
edges and nodes into the network at each step. This modified model relaxes the requirement of having 
the complete information on network evolution while doing inference and can deal with cases where we only
have coarse timestamps available.

\subsection{Traditional directed PA model}
We first give a description of  a  traditional directed PA
model, where only one new edge is created at each step; 
a formal construction of this traditional model is deferred to Section~\ref{subsec:PA_setup}.
This is also a special case of the directed PA model considered in \cite{wan:wang:davis:resnick:2017, krapivsky:redner:2001, bollobas:borgs:chayes:riordan:2003}.

We initialize the model with graph $\overline{\bfG}(0)$, which consists of one node (labeled as Node 1) and a self-loop.
Let $\overline{\bfG}(n)$ denote the graph  after $n$ steps
and $\overline{\bfV}(n)$ be the set of nodes in $\overline{\bfG}(n)$ with $\overline{\bfV}(0) = \{1\}$ and $|\overline{\bfV}(0)| = {1}$.
Set $\bigl(\barI_v(n), \barO_v(n)\bigr)$ to
be the in- and out-degrees of node $v$ in $\overline{\bfG}(n)$.

At each step, with probability $p$, we add a new edge which
starts from a new node and points to one of the existing nodes $v$, and the existing node is chosen with probability 
\beqq\label{eq:PA_in}
\frac{\barI_v(n)+\deltain}{\sum_{v\in\overline{\bfV}(n)} (\barI_v(n)+\deltain)}. 
\eeqq
With probability $1-p$, a new edge is added between two existing nodes $w\mapsto v$, where the starting and the ending nodes $w, v$ are chosen independently with probability
\beqq\label{eq:PA_beta}
\frac{\barI_v(n)+\deltain}{\sum_{v\in\overline{\bfV}(n)} (\barI_v(n)+\deltain)}\frac{\barO_w(n)+\deltaout}{\sum_{w\in\overline{\bfV}(n)} (\barO_w(n)+\deltaout)}.
\eeqq
Note that in the traditional PA set up, for $n\ge 0$, we have
\beqq\label{eq:sum_PA}
\sum_{v\in\overline{\bfV}(n)} \barI_v(n) = n+1 = \sum_{v\in\overline{\bfV}(n)}\barO_v(n),
\eeqq
Then the attachment probabilities in \eqref{eq:PA_in} and \eqref{eq:PA_beta} become
\[
\frac{\barI_v(n)+\deltain}{n+1+\deltain |\barbV(n)|}\qquad\text{and}\qquad
\frac{\barI_v(n)+\deltain}{n+1+\deltain |\barbV(n)|}\frac{\barO_v(n)+\deltaout}{n+1+\deltaout |\barbV(n)|},
\]
respectively.

The total number of nodes in $\overline{\bfG}(n)$ then
satisfies
\[
\frac{|\overline{\bfV}(n)|}{n}\convas p,\qquad n\to\infty.
\]
The asymptotic limit of {empirical frequencies}
\[
\frac{1}{|\barbV(n)|}\sum_{v\in \barbV(n)}\ind_{\left\{\bigl(\barI_v(2n), \barO_v(2n)\bigr) = (m,l)\right\}}
\]
in the traditional directed PA model has been studied in \cite{krapivsky:redner:2001, bollobas:borgs:chayes:riordan:2003, resnick:samorodnitsky:towsley:davis:willis:wan:2016, wang:resnick:2019}, and the
statistical fitting of this traditional PA model is discussed in \cite{wan:wang:davis:resnick:2017}.

In this model, edges are added sequentially so it fails to accommodate the real-data scenario where only coarse timestamp
information is observed. Therefore, we propose a modified PA model in the next section.

\subsection{PA model with Poisson measurement}\label{subsec:POPA}
We describe a modified directed PA model, which is 
a sequence of growing graphs $\{\bfG(n): n\geq 0\}$ with node set $\{\bfV(n): n \geq 0\}$ such that
the graph, $\bfG(0)$, starts with one node (labeled as Node 1) and a self-loop.
The formal construction of this model is given in Section~\ref{subsec:def_POPA}.
We use $\bfV(0)$ to denote the set of nodes in $\bfG(0)$ so that $\bfV(0) = \{1\}$ and $|\bfV(0)| = 1$.
From $\bfG(n)$ to $\bfG(n+1)$, $n\ge 0$, we assume the network keeps growing such that the 
number of newly created edges is always greater than or equal to 1, which agrees with observations from 
real datasets, e.g. Facebook and Slashdot. 
Motivated by the findings summarized in Section~\ref{subsec:intro_fb}, we assume the 
number of new edges from $\bfG(n)$ to $\bfG(n+1)$, denoted by $\Delta M_{n+1}$, follows a \emph{unit-shifted Poisson distribution} 
with pmf 
\beqq\label{eq:pmf_sPois}
\PP(\Delta M_{n+1} = k) = e^{-\lambda}\frac{\lambda^{k-1}}{(k-1)!},\qquad k\ge 1,\, \lambda>0,
\eeqq
and $\{\Delta M_n: n\ge 1\}$ are iid.
From $\bfG(n)$ to $\bfG(n+1)$, $n\ge 0$, we observe $\Delta M_{n+1}$ (independent from $\bfG(n)$) new edges which are created following a preferential attachment rule outlined below. 

Write $M_n := \sum_{k=1}^n \Delta M_k$. For all of the $\Delta M_{n+1}$ newly created edges,
there are two possibilities for how a new edge is added:
\begin{enumerate}
\item[(i)] With probability $p$, the new edge
starts from a new node and points to one of the existing nodes $v\in \bfV(n)$, where the existing node is chosen with probability 
\beqq\label{eq:POPA_in}
\frac{I_v(n)+\deltain}{\sum_{v\in \bfV(n)} (I_v(n)+\deltain)}.
\eeqq
\item[(ii)] With probability $1-p$, a new edge linking two existing nodes $w\mapsto v$  is created, where the starting and ending nodes $w, v$ are chosen independently with probability
\begin{align}
&\frac{I_v(n)+\deltain}{\sum_{v\in \bfV(n)} (I_v(n)+\deltain)}\frac{O_w(n)+\deltaout}{\sum_{w\in \bfV(n)} (O_w(n)+\deltaout)}.
\label{eq:POPA_beta}
\end{align}
\end{enumerate}
Note also that in $\bfG(n)$, $n\ge 1$,
\beqq\label{eq:sum_POPA}
\sum_{v\in \bfV(n)} I_v(n) = 1+\sum_{k=1}^n \Delta M_k = \sum_{v\in \bfV(n)} O_v(n),
\eeqq
so the attachment probabilities in \eqref{eq:POPA_in} and \eqref{eq:POPA_beta} become
\[
\frac{I_v(n)+\deltain}{1+M_n+\deltain|\bfV(n)|},\quad\text{and}\quad
\left(\frac{I_v(n)+\deltain}{1+ M_n+\deltain|\bfV(n)|}\right)\left(\frac{O_v(n)+\deltaout}{1+M_n+\deltaout|\bfV(n)|}\right),
\]
respectively.
In either case, the probability of having a new edge pointing to $v\in \bfV(n)$ is equal to
\begin{align*}
\frac{I_v(n)+\deltain}{1+M_n+\deltain|\bfV(n)|}.
\end{align*}
The attachment probabilities remain fixed until all $\Delta M_{n+1}$ edges are added.

The model setup given above assures that as $n\to\infty$,
\[
\frac{|\bfV(n)|}{n}\convas (\lambda+1)p.
\]
Note that in this new PA model, how the new $\Delta M_{n+1}$ edges are added depends on the configuration of $\bfG(n)$,
which addresses the possibility of having coarse timestamp information as pointed out in Section~\ref{sec:intro_data}.

\section{Traditional PA Model}\label{sec:embed_PA}
We start with the theoretical analysis on the traditional PA model, and then move to the PA model with Poisson measurement by 
modifying the attachment probabilities. Studies on the asymptotic properties of the traditional PA model can be found in 
\cite{bollobas:borgs:chayes:riordan:2003, resnick:samorodnitsky:towsley:davis:willis:wan:2016,wang:resnick:2019,
vanderHofstad:2017},
and issues with regard to the statistical inference on the traditional PA model are discussed in \cite{wan:wang:davis:resnick:2017, wan:wang:davis:resnick:2017b}.

First, we point out that the directed PA model studied in \cite{wang:resnick:2019} is a special case of the traditional PA model considered in \cite{wan:wang:davis:resnick:2017, krapivsky:redner:2001, bollobas:borgs:chayes:riordan:2003, resnick:samorodnitsky:towsley:davis:willis:wan:2016}.
It adds one new edge at each step, 
and this new edge either goes from the new node to one of the existing nodes or from one existing node to the new one. 
Then the degree sequence is embedded into a sequence of paired \emph{switched birth processes
with immigration} (SBI processes), from which the direction of a new edge is determined by the competition among exponential clocks. 
However, the traditional PA model summarized in Section~\ref{sec:prelim} is different as it allows new edges to be added between two existing nodes,
and the two existing nodes are chosen independently given the configuration in $\overline{G}(n)$. This makes the SBI embedding method inapplicable.
To overcome such problem in embedding, we now use a different construction described below.

\subsection{Model construction}\label{subsec:PA_setup}
On $\left(\mathbb{N}^2\right)^\infty$, we construct sequentially a paired process, which serves as the
in- and out-degree sequences in a traditional PA model, 
$$\left\{(\barbfI(n), \barbfO(n)): n\ge 0\right\} := 
\left\{\bigl(\barI_v(n), \barO_v(n)\bigr)_{v\ge 1}: n\ge 0\right\},$$
The construction uses the notation: for
  $v\geq 1,$
\begin{align*}
\mathbf{e}_v &:= \bigl(\underbrace{(0,0)}_{\text{first entry}},\ldots, (0,0),\underbrace{(1,0)}_{\text{$v$-th entry}}, (0,0),\ldots\bigr)\\
\widetilde{\mathbf{e}}_v &:= \bigl(\underbrace{(0,0)}_{\text{first entry}},\ldots, (0,0),\underbrace{(0,1)}_{\text{$v$-th entry}}, (0,0),\ldots\bigr).
\end{align*}
We also define two sequences of choice variables $\{J_n: n\ge 1\}$ and
$\{\widetilde J_n: n\ge 1\}$, tracking {which nodes have}  in- and out-degrees 
increased at each step. \sid{Specifying} the distribution of $\{J_n: n\ge 1\}$ and
$\{\widetilde J_n: n\ge 1\}$ \sid{uses the following ingredients.}
\sid{Let}  $\{\tau_{n,v},\widetilde{\tau}_{n,v}: v\ge 1, n\ge 1\}$
be iid unit rate exponential random variables  and  let $\{B_k: k\ge 1\}$ be iid Bernoulli random variables with $\PP(B_k = 0) =1-p = 1-\PP(B_k=1)$, 
which are independent from $\{\tau_{n+1,v},\widetilde{\tau}_{n+1,v}:
v\ge 1, n\ge 0\}$.

\sid{As in Section~\ref{sec:prelim}, assume the initialization  as 
\[
\bigl(\barbfI(0), \barbfO(0)\bigr) 
= \bigl((1,1), (0,0),(0,0),\ldots\bigr),
\]
which corresponds to a single initial node with a self-loop.}
Write $\barbV(0) := \{1\}$.
For $n=1$, set $J_1 = 1$, and
\[
\bigl(\barbfI(1), \barbfO(1)\bigr) 
= \bigl(\barbfI(0), \barbfO(0)\bigr) + \mathbf{e}_{J_1},
\]
which corresponds to increasing the in-degree of Node 1 by 1.
Then for $n=2$, let $\widetilde J_1 := 1$, $\barbV(1) := \{1, 1+B_1\}$, and 
\[
\bigl(\barbfI(2), \barbfO(2)\bigr) 
= \bigl(\barbfI(1), \barbfO(1)\bigr) + (1-B_1)\widetilde{\mathbf{e}}_{\widetilde J_1} + B_1\widetilde{\mathbf{e}}_2,
\]
which corresponds to incrementing the out-degree of Node 1 by 1 with probability $1-p$ and adding a node with out-degree 1
with probability $p$.

For $n\ge 1$, 
we set $\barbV(n) := \{1,2,\ldots,1+\sum_{k=1}^n B_k\}$, and
define
$J_{n+1}$ and $\widetilde J_{n+1}$ as 
\begin{align}
J_{n+1} &:= \argmin_{v\in \barbV(n)} \frac{\tau_{n+1,v}}{\barI_v(2n)+\deltain},\label{eq:def_J}\\
\widetilde{J}_{n+1} &:= \argmin_{v\in \barbV(n)} \frac{\widetilde\tau_{n+1,v}}{\barO_v(2n)+\deltaout}.\label{eq:def_tildeJ}
\end{align}
So for each fixed $n$, $J_{n+1}$ and $\widetilde J_{n+1}$ are conditionally independent given 
the filtration
$\mathcal{F}_{2n} := \sigma\left\{(\barbfI(k), \barbfO(k)): 0\le k\le 2n\right\}$.
Since for $n\ge 1$, $1+\sum_{k=1}^n B_k$ is equal to the number of 
$\barO_v(\cdot)$, $v\in \barbV(n)$, that has non-zero values, 
 $\{B_k: 1\le k\le n\}$ is $\mathcal{F}_{2n}$-measurable.
Having defined $\{\bigl(\barbfI({k}), \barbfO({k})\bigr): 1\le k\le 2n\}$, \sid{set}
\begin{align}
\bigl(\barbfI({2n+1}),& \barbfO({2n+1})\bigr) 
= \bigl(\barbfI({2n}), \barbfO({2n})\bigr) + \mathbf{e}_{J_{n+1}},\label{eq:def_barI}\\
  \intertext{which corresponds to incrementing the in-degree of
   existing node $J_{n+1}$, and}
\bigl(\barbfI({2n+2}),& \barbfO({2n+2})\bigr) \nonumber\\
&= \bigl(\barbfI({2n+1}), \barbfO({2n+1})\bigr) + (1-B_{n+1})\widetilde{\mathbf{e}}_{\widetilde J_{n+1}} 
+ B_{n+1}\widetilde{\mathbf{e}}_{1+|\barbV(n)|},\label{eq:def_barO}
\end{align}
which corresponds to incrementing the out-degree of an existing node according to the choice variable $\widetilde J_{n+1}$
with probability $1-p$ and adding a node with out-degree 1 with probability $p$.
From this construction, we have for $n\ge 0$,
\[
\sum_{v\in \barbV(n)} \barI_v(2n) = n+1 = \sum_{v\in \barbV(n)} \barO_v(2n),
\]
which agrees with \eqref{eq:sum_PA}.

\sid{For fixed $n$}, we consider
$\bigl(\barI_v(2n), \barO_v(2n)\bigr)$, $v\in \barbV(n)$, as the in-
and out-degrees of Node $v$ in $\overline{\bfG}(n)$, and 
 $\bigl(\barbfI({2n}), \barbfO({2n})\bigr)$ is the in- and
out-degree sequences in the traditional PA model. 
We write $\PP^{\mathcal{F}_{2n}}(\cdot):= \PP(\cdot|\mathcal{F}_{2n})$, then
for $n\ge 0$, the transition probability from 
$\bigl(\barbfI(2n), \barbfO(2n)\bigr)$ to $\bigl(\barbfI(2(n+1)), \barbfO(2(n+1))\bigr)$ becomes:
for $v,w\in \barbV(n)$,
\begin{align}
\PP^{\mathcal{F}_{2n}}&\left(\bigl(\barbfI(2(n+1)), \barbfO(2(n+1))\bigr) = \bigl(\mathbf{I}(2n), \barbfO(2n)
+ \mathbf{e}_v +  \widetilde{\mathbf{e}}_{1+|\barbV(n)|}\bigr)\right) \nonumber\\
&= \PP^{\mathcal{F}_{2n}}\left(J_{n+1}=v, B_{n+1}= 1\right)\nonumber\\
&= p\, \frac{\barI_v({2n})+\deltain}{n+1+\deltain |\barbV(n)|},\label{eq:J_PA}
\intertext{which agrees with the scenario described in \eqref{eq:PA_in}, and}
\PP^{\mathcal{F}_{2n}}&\left(\bigl(\barbfI(2(n+1)), \barbfO(2(n+1))\bigr) = \bigl(\barbfI(2n), \barbfO(2n)
+ \mathbf{e}_v +  \widetilde{\mathbf{e}}_{w}\bigr)\right)\nonumber \\
&= \PP^{\mathcal{F}_{2n}}\left(J_{n+1}=v, \widetilde{J}_{n+1}=w, B_{n+1}= 0\right)\nonumber\\
&= (1-p)\, \left(\frac{\barI_v({2n})+\deltain}{n+1+\deltain |\barbV(n)|}\right)\left(\frac{\barO_w({2n})+\deltaout}{n+1+\deltaout |\barbV(n)|}\right),\label{eq:Jbeta_PA}
\end{align}
which agrees with the second scenario in \eqref{eq:PA_beta}.
Equations \eqref{eq:J_PA} and \eqref{eq:Jbeta_PA} also show that for $n\ge 1$, both
$J_n$ and $\widetilde{J}_n$ are independent from $B_n$.

By the definition of $\{J_k: k\ge1\}$ and $\{\widetilde{J}_k: k\ge1\}$, we see that
for \sid{$v_{i+1}, w_{i+1}\in \barbV(i)$, $i\ge 0$, }
\[
\PP^{\mathcal{F}_{2n}}\left(J_{n+1} = v_{n+1}, \widetilde{J}_{n+1}=w_{n+1}\right)
= \PP^{\mathcal{F}_{2n}}\left(J_{n+1} = v_{n+1}\right)\PP^{\mathcal{F}_{2n}}\left(\widetilde{J}_{n+1}=w_{n+1}\right).
\]
Then we have for $v_n, w_n \in \barbV(n-1)$, $n\ge 1$, 
\begin{align}
&\PP^{\mathcal{F}_{2(n-1)}\vee B_n}\left(J_{n} = v_{n}, J_{n+1} =
                v_{n+1},\widetilde{J}_{n}=w_n,
                \widetilde{J}_{n+1}=w_{n+1}\right) \nonumber \\
&= \EE^{\mathcal{F}_{2(n-1)}\vee B_n}\left(\ind_{\{J_{n} = v_{n}\}}\PP^{\mathcal{F}_{2n}}\left(J_{n+1} = v_{n+1}\right)
\ind_{\{\widetilde J_{n} = w_{n}\}}\PP^{\mathcal{F}_{2n}}\left(\widetilde J_{n+1} = w_{n+1}\right)
\right).\nonumber \\
\intertext{By \eqref{eq:def_J} and \eqref{eq:def_tildeJ}, we have
$\barI_{v_{n+1}}(2n) = \barI_{v_{n+1}}(2(n-1))+\ind_{\{J_{n} = v_{n+1}\}}$ and 
$\barO_{w_{n+1}}(2n) = \barO_{w_{n+1}}(2(n-1))+\ind_{\{\widetilde J_{n} = w_{n+1}\}}$, so that}
&\PP^{\mathcal{F}_{2(n-1)}\vee B_n}\left(J_{n} = v_{n}, J_{n+1} =
                                                                                                  v_{n+1},\widetilde{J}_{n}=w_n, \widetilde{J}_{n+1}=w_{n+1}\right)\nonumber \\ 
&=  \EE^{\mathcal{F}_{2(n-1)}\vee B_n}\left(
\ind_{\{J_{n} = v_{n}\}} \frac{\barI_{v_{n+1}}(2(n-1))+\deltain +
                                                                                                                                                                                   \ind_{\{J_{n}
                                                                                                                                                                                   =
                                                                                                                                                                                   v_{n+1}\}}}{n+1+\deltain
                                                                                                                                                                                   |\barV(n)|}\right. \nonumber
  \\
&\left.\qquad\quad\times
\ind_{\{\widetilde J_{n} = w_{n}\}}
       \frac{\barO_{w_{n+1}}(2(n-1))+\deltaout + \ind_{\{\widetilde
       J_{n} = w_{n+1}\}}}{n+1+\deltaout |\barV(n)|} 
    \right)\nonumber \\
  \intertext{\sid{
  and since $\ind_{\{ J_n=v_n\}}
  \ind_{\{  J_{n}=v_{n+1}\}} = \ind_{\{ J_n=v_n\}}
  \ind_{\{v_n=v_{n+1}  \}}$ with a similar result for
  $\widetilde{J}_n$ this is,}}
&= \frac{\barI_{v_{n+1}}(2(n-1))+\deltain + \ind_{\{v_{n} = v_{n+1}\}}}{n+1+\deltain |\barV(n)|}
\frac{\barO_{w_{n+1}}(2(n-1))+\deltaout + \ind_{\{w_{n} =
                                 w_{n+1}\}}}{n+1+\deltaout
                                 |\barV(n)|}\nonumber \\ 
&\qquad\times \PP^{\mathcal{F}_{2(n-1)}}\left(J_{n} =
                                                            v_{n}\right)\PP^{\mathcal{F}_{2(n-1)}}\left(\widetilde{J}_{n}=w_{n}\right), \label{e:jt} 
\end{align}
\sid{where the last step used the independence of $B_n$ from $(J_n,
  \widetilde{J}_n)$. Taking marginals in \eqref{e:jt} gives}
\begin{align*}
\PP^{\mathcal{F}_{2(n-1)}\vee B_n}&\left(J_{n} = v_{n}, J_{n+1} = v_{n+1}\right)\\
 =& \frac{\barI_{v_{n+1}}(2(n-1))+\deltain + \ind_{\{v_{n} = v_{n+1}\}}}{n+1+\deltain |\barV(n)|}
 \frac{\barI_{v_{n}}(2(n-1))+\deltain }{n+1+\deltain |\barV(n-1)|}\\
 =& \frac{\barI_{v_{n+1}}(2(n-1))+\deltain + \ind_{\{v_{n} = v_{n+1}\}}}{n+1+\deltain |\barV(n)|}
 \PP^{\mathcal{F}_{2(n-1)}}\left(J_{n} = v_{n}\right),
\end{align*}
and similarly,
\begin{align*}
\PP^{\mathcal{F}_{2(n-1)}\vee B_n}&\left(\widetilde J_{n} = w_{n}, \widetilde J_{n+1} = w_{n+1}\right)\\
=& \frac{\barO_{w_{n+1}}(2(n-1))+\deltaout + \ind_{\{w_{n} = w_{n+1}\}}}{n+1+\deltaout |\barV(n)|}
\PP^{\mathcal{F}_{2(n-1)}}\left(\widetilde{J}_{n}=w_{n}\right).
\end{align*}
\sid{Therefore,}
\begin{align*}
&\PP^{\mathcal{F}_{2(n-1)}\vee B_n}\left(J_{n} = v_{n}, J_{n+1} = v_{n+1},\widetilde{J}_{n}=w_n, \widetilde{J}_{n+1}=w_{n+1}\right)\\
&= \PP^{\mathcal{F}_{2(n-1)}\vee B_n}\left(J_{n} = v_{n}, J_{n+1} = v_{n+1}\right)
\PP^{\mathcal{F}_{2(n-1)}\vee B_n}\left(\widetilde{J}_{n}=w_n, \widetilde{J}_{n+1}=w_{n+1}\right).
\end{align*}
Hence, following an induction argument in $n$, 
we see that for $v_k, w_k \in \barbV(k-1)$, $1\le k\le n$,
\begin{align}
\label{eq:Js_indep}
\PP^{\{B_k\}_{k=1}^{n}}&\left(J_{k} = v_{k},\widetilde{J}_{k}=w_k, 1\le k\le n\right)\nonumber\\
= &\PP^{\{B_k\}_{k=1}^{n}}\left(J_{k} = v_{k}, 1\le k\le n\right)\PP^{\{B_k\}_{k=1}^{n}}\left(\widetilde{J}_{k}=w_k, 1\le k\le n\right).
\end{align}

\subsection{Degree Distribution}\label{subsec:PA_dist}
In this section, we study the in- and out-degree distribution in a
traditional PA model by embedding each of them into a sequence of
birth-immigration processes. This will serve as a model for how to
analyze the degree distribution {of} a PA model with Poisson
measurement. 

From the construction in the previous section, we see that
with the additional definition $J_0 := 1$ and $\widetilde J_0:=1$,
\begin{align*}
\left\{\bigl(\barI_v(2n), \barO_v(2n)\bigr): v\in \barbV(n), n\ge 0\right\} = \left\{\left(\sum_{k=1}^n \ind_{\{J_k = v\}}, \sum_{k=1}^n \ind_{\{\widetilde J_k = v\}}\right): v\in \barbV(n), n\ge 0\right\}.
\end{align*}
By \eqref{eq:Js_indep}, $\{J_k: 1\le k\le n\}$ and $\{\widetilde J_k: 1\le k\le n\}$ are conditionally independent given $\{B_k: 1\le k\le n\}$.
We first
specify the marginal distributions of $\{\barbfI(2n): n\ge 0\}$ and $\{\barbfO(2n): n\ge 0\}$, conditional on $\{B_k: k\ge 1\}$.

Let $S_v$, $v\in \barbV(n)$, be the number of edges that have been added into the graph when the $v$-th node is created, i.e.
$S_1 = 0$, and
\beqq\label{eq:def_Sv}
S_v:=\inf\left\{n\ge 1: 1+\sum_{k=1}^n B_k = v\right\},\qquad v\ge 2.
\eeqq
Equation~\eqref{eq:def_Sv} reveals that $S_v$ is the waiting time in Bernoulli trials until $v-1$ successes have been achieved.
Hence, $S_v$ follows a negative binomial distribution with generating function
\[
\EE(s^{S_v})= (s+(1-s)/p)^{-(v-1)},\qquad s\in [0,1].
\]
We then have
\[
\left\{\bigl(\barI_v(2n), \barO_v(2n)\bigr): v\in \barbV(n), n\ge 0\right\} = \left\{\left(\sum_{k=S_v}^n \ind_{\{J_k = v\}}, \sum_{k=S_v}^n \ind_{\{\widetilde J_k = v\}}\right): v\in \barbV(n), n\ge 0\right\}.
\]

\subsubsection{Embedding}\label{subsubsec:embed}

The key ingredient used to specify the marginal distributions is the  framework built from birth-immigration processes (cf. \cite{wang:resnick:2018, athreya:ghosh:sethuraman:2008}), where the choice variables, $\{J_k: k\ge 1\}$ and $\{\widetilde J_k: k\ge 1\}$, can be viewed as marking which
birth immigration process jumps first.
We now start with a brief overview on the birth immigration process.
A linear birth-immigration process, 
$\{BI_\delta (t): t\ge 0\}$, having unit lifetime parameter and
immigration parameter $\delta\ge 0$ is a continuous time Markov
process with state space $\mathbb{N}$ and
transition rate  
$$q_{k,k+1} = k + \delta,\qquad k\ge 0.$$
When $\delta = 0$ there is no immigration and the birth-immigration process becomes
a pure birth process and in such cases, the process usually starts from 1.
For $\delta>0$, the birth-immigration process starting from 0 population can be
constructed from a Poisson process and an independent family of iid
linear birth processes \cite{tavare:1987}.

{Following the procedure in \cite{wang:resnick:2018},}
we embed $\{\barbfI(2n): n\ge 0\}$ into a sequence of birth immigration processes, which are 
independent from the Bernoulli random variables $\{B_k: k\ge 1\}$.
Suppose that $\{BI^{(v)}_{\ind_{\{v=1\}}+\deltain}(t): t\ge 0\}_{v\ge 1}$ is a
sequence of independent birth-immigration processes, all of which start with population 0, have a unit lifetime parameter and immigration parameters equal to $1+\deltain$ for $BI^{(1)}_{1+\deltain}$ and $\deltain$ for $BI^{(v)}_{\deltain}$, $v\ge 2$. 
At time $\Gamma_0:=0$, we start with having only $\{BI^{(1)}_{1+\deltain}(t): t\ge 0\}$, 
and let $\Gamma_1$ be the time at which
the first jump of $BI^{(1)}_{1+\deltain}(\cdot)$ occurs. Set also $J'_1 =1$, representing that $BI^{(1)}_{1+\deltain}(\cdot)$ jumps to 1 at $\Gamma_1$.

At time $\Gamma_1$, 
if $B_1 = 1$ start a new birth-immigration process
$\{BI^{(1+B_1)}_{\deltain}(t\sid{-\Gamma_1}): t\ge \Gamma_1\}$,
and use $\Gamma_2$ to denote the first time after $\Gamma_1$ such that one of the $BI^{(1)}_{1+\deltain}(\cdot)$ and $BI^{(1+B_1)}_{\deltain}(\cdot)$
processes jumps. 
If $B_1=0$, then $\Gamma_2$ is the first time after $\Gamma_1$ such that the $BI^{(1)}_{1+\deltain}(\cdot)$ process jumps,
and no new birth-immigration process is initiated.
Write $\barbV(1) := \{1, 1+B_1\}$, \sid{so $|\barbV(1)|=1+B_1$ is the number
of processes running at $\Gamma_1$.}  Let $J'_2$ denote which birth-immigration process  jumps at $\Gamma_2$.
For $n\ge 1$,  proceed in the way \sid{just} outlined:   Given $\{B_k: k=1,\ldots, n\}$,
$\barbV(n) := \{1, \ldots, 1+\sum_{k=1}^nB_k\}$ \sid{indexes the processes running,}
 $\Gamma_{n+1}$ is the first time after $\Gamma_n$ when one
of the processes
$$\{BI^{(v)}_{\ind_{\{v=1\}}+\deltain}(t-\Gamma_{S_v}): t\ge \Gamma_{S_v}\},
v\in\barbV(n),$$ 
jumps, 
and $J_{n+1}'$  denotes which process jumps at $\Gamma_{n+1}$.
Then start a new birth-immigration process $\{BI^{(|\barbV(n+1)|)}_{\deltain}(t-\Gamma_{n+1}): t\ge \Gamma_{n+1}\}$ at $\Gamma_{n+1}$ if $B_{n+1} = 1$.
By \cite[Proposition~2.1]{athreya:ghosh:sethuraman:2008}, we see that given $\{B_k: 1\le k\le n\}$,
$\{\Gamma_{k+1}-\Gamma_k: 0\le k\le n\}$ are independent exponential random variables with means
$(k+1+\deltain |\barbV(k)|)^{-1}$, $0\le k\le n$.

Similar to the embedding results in \cite{wang:resnick:2018}, we have, for
\begin{equation}\label{e:sigmafield}
  \mathcal{F}'_n :=
  \sigma\left(\left\{B_k\right\}_{k=1}^n;\left\{BI^{(v)}_{\ind_{\{v=1\}}+\deltain}(t-\Gamma_{S_v}):
      \Gamma_{S_v}\le t\le \Gamma_n\right\}_{v\in\barbV(n)}\right),\end{equation}
 and $v_k\in\barbV(n)$, $t_k\ge 0$, $k=1,\ldots, n+1$,
\begin{align*}
&\PP^{\{B_k\}_{k=1}^{n+1}}\left(\bigcap_{k=1}^{n+1} \left\{J'_{k}=v_k, \Gamma_{k}-\Gamma_{k-1}> t_{k}\right\}\right)\\
&=\EE^{\{B_k\}_{k=1}^{n+1}} \left(\ind_{\left\{\bigcap_{k=1}^{n} \left\{J'_{k}=v_k, \Gamma_{k}-\Gamma_{k-1}> t_{k}\right\}\right\}}
\PP^{\mathcal{F}'_n}\left(J'_{n+1}=v_{n+1}, \Gamma_{n+1}-\Gamma_n> t_{n+1}\right)\right)\\
&= \EE^{\{B_k\}_{k=1}^{n+1}} \left(\ind_{\left\{\bigcap_{k=1}^{n} \left\{J'_{k}=v_{n+1}, \Gamma_{k}-\Gamma_{k-1}> t_{k}\right\}\right\}}
\PP^{\mathcal{F}'_n}\left(J'_{n+1}=v_{n+1}\right)\PP^{\mathcal{F}'_n}\left(\Gamma_{n+1}-\Gamma_n> t_{n+1}\right)\right).
\end{align*}
So conditional on $\{B_k: k\ge 1\}$, we have
\begin{align}
\left\{(J'_k,\Gamma_{k}-\Gamma_{k-1}) : k\ge 1\right\}&\stackrel{d}{=} \left\{\left(J_k,\bigwedge_{v\in \barbV(k-1)} \frac{\tau_{k,v}}{\barI_v(2(k-1))+\deltain}\right): k\ge 1\right\}.\label{eq:JGam_PA}
\end{align}
Now applying the embedding framework in \cite[Theorem 3]{wang:resnick:2018} gives,
\beqq\label{eq:in_dist_PA}
\{\barbfI(2n): n\ge 0\} \stackrel{d}{=} \left\{\left\{BI^{(v)}_{\ind_{\{v=1\}}+\deltain}(\Gamma_n-\Gamma_{S_v}): v\in \barbV(n)\right\}: n\ge 0\right\}.
\eeqq

The \sid{embedding} of $\{\barbfO(2n): n\ge 0\}$ follows in a similar way.
First, note that for all $v\in \barbV(n)$, $n\ge 0$, $\barO_v(2n)\ge 1$, but similar to the in-degree case, we here assume the birth-immigration process starts with population 0. So we now actually embed $\left\{\barbfO(2n)-\mathbbm{1}_{|\barbV(n)|}: n\ge 0\right\}$ into a sequence of birth-immigration process, where
$
\mathbbm{1}_n := (\underbrace{1,1,\ldots,1}_{\text{n entries}},0,\ldots).
$
Suppose that $\{\widetilde{BI}^{(v)}_{1+\deltaout}(t): t\ge 0\}_{v\ge 1}$ is a
sequence of independent birth-immigration processes, all of which
are independent from $\{{BI}^{(v)}_{\ind_{\{v=1\}}+\sid{\deltain}}(\cdot)\}_{v\ge 1}$
and $\{B_k:k\ge 1\}$,
 start with population 0, have a unit lifetime parameter and immigration parameters equal to $1+\deltaout$. 
At time $\widetilde{\Gamma}_0:=0$, \sid{initiate} $\{\widetilde{BI}^{(1)}_{1+\deltaout}(t): t\ge 0\}$. 
If $B_1= 1$, start another process $\{\widetilde{BI}^{(1+B_1)}_{1+\deltaout}(t): t\ge 0\}$ at time 0,
which corresponds to
adding a new node born with out-degree 1 with probability $p$. 
Here although we do not have Node 0 in the PA model, we still set 
$\widetilde J_1' = 0$ to represent the situation where none of the existing nodes has changes in their out degrees.
Otherwise, if $B_1=0$, 
let $\widetilde T_1$ be the time at which
the first jump of $\widetilde{BI}^{(1)}_{1+\deltaout}(\cdot)$ occurs
and set $\widetilde{J}'_1 =1$. 
This corresponds to incrementing the out-degree of Node 1 by 1 with probability $1-p$, 
but no new node is added.
In order to combine two scenarios, we further set $\widetilde T_0 := 0$ and
define $\widetilde \Gamma_1 := \widetilde T_{1-B_1}$,
then equivalently, we initiate a new birth-immigration process
$\{\widetilde{BI}^{(1+B_1)}_{1+\deltaout}(t-\sid{\widetilde\Gamma_1)}: t\ge \widetilde\Gamma_1\}$,
at time $\widetilde\Gamma_1$.

Next, if $B_2=0$, let $\widetilde T_2$ denote the first time after $\widetilde\Gamma_1$ 
such that one of the $\widetilde{BI}^{(1)}_{1+\deltaout}(\cdot)$ and $\widetilde{BI}^{(1+B_1)}_{1+\deltaout}(\cdot-\widetilde \Gamma_1)$
processes jumps,
then
\[
\PP\left(\widetilde{T}_2-\widetilde\Gamma_1>t\middle|B_2=0, B_1\right)
= e^{-(2+\deltaout(1+B_1))t},\qquad t\ge 0.
\]
Let $\widetilde{J}'_2$ denote which birth-immigration process jumps at $\widetilde T_2$.
If instead $B_2=1$, we write $\widetilde{J}_2'=0$, and initiate a new birth-immigration process $\{\widetilde{BI}^{\left(1+\sum_{k=1}^2 B_k\right)}_{1+\deltaout}(t-{\widetilde\Gamma_1)}: t\ge \widetilde\Gamma_1\}$ at time $\widetilde\Gamma_1$. 
Using the consolidating notation
 $\widetilde\Gamma_2 := \widetilde{T}_{\sum_{k=1}^2(1-B_k)} = \widetilde{T}_{3-|\barbV(2)|}$,
the two-scenario procedure described above is equivalent to
 initiate a new birth-immigration process 
$\{\widetilde{BI}^{\left(1+\sum_{k=1}^2 B_k\right)}_{1+\deltaout}(t-{\widetilde\Gamma_2)}: t\ge \widetilde\Gamma_2\}$
at $\widetilde\Gamma_2$.

For $n\ge 1$, we proceed in the way outlined above such that given $\{B_k: k=1,\ldots, n\}$ and
\beqq\label{eq:BI_PAout}
\left\{\widetilde{BI}^{(v)}_{1+\deltaout}(t-\widetilde\Gamma_{S_v}): t\ge \widetilde\Gamma_{S_v}
  \right\}, v\in \barbV(n),
\eeqq
if $B_{n+1}=0$, we set $\widetilde T_{n+2-|\barbV(n)|}$ to be the first time after $\widetilde\Gamma_n
\equiv \widetilde T_{n+1-|\barbV(n)|}$ when one of the processes in \eqref{eq:BI_PAout} jumps, and
use $\widetilde{J}_{n+1}'$ to denote which process in \eqref{eq:BI_PAout} jumps at $\widetilde T_{n+2-|\barbV(n)|}$.
If $B_{n+1} = 1$,
we set $\widetilde{J}_{n+1}' = 0$, and
 start a new birth-immigration process $\{\widetilde{BI}^{(|\barbV(n)|+1)}_{1+\deltaout}(t-\widetilde\Gamma_{n}): t\ge \widetilde\Gamma_{n}\}$ at $\widetilde\Gamma_{n}$.
 Also, write $\widetilde\Gamma_{n+1} := \widetilde T_{n+2-|\barbV(n+1)|}$,
and we have that conditional on $\{B_k: k\ge 1\}$, 
\begin{align*}
\left\{(\widetilde J'_k, \,\widetilde\Gamma_{k}-\widetilde\Gamma_{k-1}): k\ge 1\right\}&\stackrel{d}{=} \left\{\left((1-B_k)\widetilde J_k,\bigwedge_{v\in \barbV(k-1)} \frac{(1-B_k)\widetilde\tau_{k,v}}{\barO_v(2(k-1))+\deltaout}\right): k\ge 1\right\},
\end{align*}
and 
\beqq\label{eq:out_dist_PA}
\{\barbfO(2n): n\ge 0\} \stackrel{d}{=} \left\{\left\{1+\widetilde{BI}^{(v)}_{1+\deltaout}(\widetilde\Gamma_n-\widetilde\Gamma_{S_v})\right\}_{ v \in \barbV(n)}: n\ge 0\right\}.
\eeqq

Combining \eqref{eq:in_dist_PA} and \eqref{eq:out_dist_PA} leads to the embedding results in
the following theorem.
\bthe\label{thm:PA_BI}
Given $\{B_k: k\ge 1\}$ and with $S_v$, $v\in \barbV(n)$ defined in \eqref{eq:def_Sv}, then on $\mathbb{N}^\infty$,
\begin{align*}
\left\{\barbfI({2n}): n\ge 0\right\}
\stackrel{d}{=}& \left\{\left(\left\{\ind_{\{v=1\}} + BI^{(v)}_{\ind_{\{v=1\}}+\deltain}\left(\Gamma_n - \Gamma_{S_v}\right)\right\}_{v\in \barbV(n)}, 0,\ldots\right): n\ge 0\right\}, \nonumber\\
\left\{\barbfO({2n}): n\ge 0\right\}
\stackrel{d}{=}& 
\left\{\left(\left\{1+\widetilde{BI}^{(v)}_{1+\deltaout}\left(\widetilde\Gamma_n - \widetilde\Gamma_{S_v}\right)\right\}_{v\in \barbV(n)}, 0,\ldots\right): n\ge 0\right\}.
\end{align*}
Also, for fixed $n\ge 0$, \sid{in $\mathbb{N}\times \mathbb{N}_{>0}$,}
\begin{align}
\left(\barbfI({2n}),\right.&\left.\barbfO(2n)\right)
\stackrel{d}{=}\nonumber\\
&\left(\left(\ind_{\{v=1\}} + BI^{(v)}_{\ind_{\{v=1\}}+\deltain}\left(\Gamma_n - \Gamma_{S_v}\right),
1+\widetilde{BI}^{(v)}_{1+\deltaout}\left(\widetilde\Gamma_n - \widetilde\Gamma_{S_v}\right)
\right)_{v\in \barbV(n)}, (0,0),\ldots\right).
\label{eq:PA_BI_dist}
\end{align}
\ethe

\subsubsection{Joint degree counts}
Using the embedding results in Theorem~\ref{thm:PA_BI}, we give the
convergence of joint in- and out-degree counts 
in a traditional PA model.
\bco\label{cor:PA_count}
Let $Z_{\delta_1}(p_1)$ and $\widetilde Z_{\delta_2}(p_2)$ be two independent negative binomial random variables with parameters $\delta_i> 0$,
$p_i\in (0, 1)$, $i=1,2$,
and generating functions 
\[
(s+(1-s)/p_i)^{-\delta_i}, \quad i=1,2,\qquad s\in [0,1].
\]
\sid{In a traditional PA model, as $n\to\infty$  we
have for $(m,l) \in \mathbb{N}\times \mathbb{N}_{>0}$}
\begin{align}
\frac{1}{|\barbV(n)|}&\sum_{v\in \barbV(n)}\ind_{\left\{\bigl(\barI_v(2n), \barO_v(2n)\bigr) = (m,l)\right\}}\nonumber\\
&\convp  \int_0^1 \PP\left(\left(Z_{\deltain}\left(t^{1/(1+\deltain p)}\right), 1+ \widetilde{Z}_{1+\deltaout}\left(t^{(1-p)/(1+\deltaout p)}\right)\right) = (m,l)\right)\dd t
=: p_{m,l}.
\label{eq:PA_pij}
\end{align}
\eco
Note that the limiting $p_{m,l}$ in \eqref{eq:PA_pij} agrees with the results in \cite{bollobas:borgs:chayes:riordan:2003,resnick:samorodnitsky:towsley:davis:willis:wan:2016}.
In Section~\ref{subsec:pf1}, we give a proof of Corollary~\ref{cor:PA_count} using the embedding results in Theorem~\ref{thm:PA_BI}, which is different from what is given in 
\cite{bollobas:borgs:chayes:riordan:2003,resnick:samorodnitsky:towsley:davis:willis:wan:2016}.
Such proof machinery is important when
we show in the next section that the right hand side of \eqref{eq:PA_pij} is also the limiting joint distribution for 
a PA model with Poisson measurement.

\section{The PA Model with Poisson Measurement}\label{sec:POPA}
Suppose now we observe a sequence of iid unit-shifted Poisson random variables with rate $\lambda>0$, $\{\Delta M_k: k\ge 1\}$,
whose pmf is given in \eqref{eq:pmf_sPois}.
We assume that $\{\Delta M_k: k\ge 1\}$ are 
 independent from $\{B_k: k\ge 1\}$. 
Write $M_0 := 0$ and $M_n = \sum_{k=1}^n \Delta M_k$, $n\ge 1$, 
and 
we add $\Delta M_n$ new edges from $\bfG(n-1)$ to $\bfG(n)$.

\sid{In a traditional PA model, attachment probabilities change as
  each new edge is added to the network.
However, in  a PA model with Poisson measurement,
attachment probabilities remain constant through the process of adding
$\Delta M_n$ edges for each $n$.}
We start with the formal construction of the PA model with Poisson measurement
using discrete indexing to describe addition of edges to a graph starting from a single node with self loop.

\subsection{Model construction}\label{subsec:def_POPA}

Similar to $(\barbfI(\cdot), \barbfO(\cdot))$, we define iteratively on $(\mathbb{N}^2)^\infty$,
$$\{(\bfI(n), \bfO(n)): n\ge 0\} := \{(I_v(n), O_v(n))_{v\ge 1}: n\ge 0\},$$ 
which serves as the in- and out-degree sequence 
in a PA model with Poisson measurement.
First, we set
\[
\bigl(\bfI(0), \bfO(0)\bigr) := \bigl((1,1), (0,0),\ldots\bigr),
\]
which corresponds to an initial node with a self loop.
Having observed $M_1 \equiv \Delta M_1$, we let 
$\mathcal{J}_{k} = 1$ and $\widetilde{\mathcal{J}}_{k}= 1$, $1\le k \le M_1$, be choice variables.
Define $\bfV(0) := \{1\}$, and
\[
\bigl(\bfI(1), \bfO(1)\bigr) = \bigl(\bfI(0), \bfO(0)\bigr) + M_1 \mathbf{e}_1
+ \sum_{k=1}^{M_1} (1-B_k) \widetilde{\mathbf{e}}_1 + \sum_{k=1}^{M_1} B_k \widetilde{\mathbf{e}}_{1+\sum_{l=1}^{k} B_l}.
\]
Here $M_1 \mathbf{e}_1$ corresponds to incrementing the in-degree of Node 1 by $M_1$, and for $1\le k\le M_1$,
$(1-B_k) \widetilde{\mathbf{e}}_1$ corresponds to incrementing the out-degree of Node 1 by $1$ with probability $1-p$,
 and $B_k \widetilde{\mathbf{e}}_{1+\sum_{l=1}^{k} B_l}$
corresponds to adding a new node with out-degree 1 with probability $p$.
This defines the graphs $\bfG(0)$ and
$\bfG(1)=(\bfV(1),E(1))$,
where
$\bfV(1)=\{1,\dots,1+ \sum_{i=1}^{M_1} B_i\}$, and 
$E(1)=\{(v_1,v_2): v_i \in \bfV(1)\}$.

For $n\ge 1$, 
assume we have defined $\bfG(0),\bfG(1),\dots, \bfG(n)$
and observed $\{\Delta M_{k}:1\le k\le n\}$ as well as the set of nodes
$\bfV(n) := \{1,\ldots, 1+\sum_{k=1}^{M_n} B_k\}$
\sid{
  and  $\{(\bfI (k), \bfO(k)),1\leq k \leq n\}$.}
We now construct $\bfG(n+1)$.
As in Section~\ref{subsec:PA_setup}  \sid{assume} $\{\tau_{n,v},\widetilde{\tau}_{n,v}: v\ge 1, n\ge 1\}$ are iid 
exponential random variables with unit rate which are independent from $\{M_k, B_k: k\ge 1\}$.
Then let $\{\mathcal J_{k}, \widetilde{\mathcal J}_{k}: M_n +1 \le k \le M_{n+1}, n\ge 0\}$ be 
choice variables tracking the nodes with which each new edge is associated
such that
for $M_n+1\le k\le M_{n+1}$, 
\begin{align*}
\mathcal{J}_{k} := \argmin_{v\in \bfV(n)}\frac{\tau_{k,v}}{I_v(n)+\deltain},\qquad &
\widetilde{\mathcal{J}}_{k} := \argmin_{v\in \bfV(n)}\frac{\widetilde\tau_{k,v}}{O_v(n)+\deltaout}.
\end{align*}
Next, we set
\begin{align*}
\bigl(\bfI(n+1), \bfO(n+1)\bigr) 
=& \bigl(\bfI(n), \bfO(n)\bigr) + \sum_{k=M_n+1}^{M_{n+1}}{\mathbf{e}}_{{\mathcal{J}}_k} \\
&+ \sum_{k=M_n+1}^{M_{n+1}}(1-B_{k})\widetilde{\mathbf{e}}_{\widetilde{\mathcal{J}}_k} 
+ \sum_{k=M_n+1}^{M_{n+1}}B_{k}\widetilde{\mathbf{e}}_{1+\sum_{l=1}^k B_l},
\end{align*}
where for each $M_n+1\le k \le M_{n+1}$,  ${\mathbf{e}}_{{\mathcal{J}}_k}$ corresponds to incrementing the
in-degree of an existing node by 1 according to ${\mathcal{J}}_k$, 
$(1-B_{k})\widetilde{\mathbf{e}}_{\widetilde{\mathcal{J}}_k} $ corresponds to
incrementing the out-degree of an existing node by 1 according to ${\widetilde{\mathcal{J}}_k}$ with probability $1-p$,
and $B_{k}\widetilde{\mathbf{e}}_{1+\sum_{l=1}^k B_l}$ corresponds to adding 
a new node with out-degree 1 with probability $p$.
\sid{It is helpful to observe for the last term that
$$
\sum_{k=M_n+1}^{M_{n+1}}B_{k}\widetilde{\mathbf{e}}_{1+\sum_{l=1}^k
  B_l}
=\sum_{j=1}^{\Delta M_{n+1}} B_{M_n +j }\widetilde{\mathbf{e}}_{|\bfV
  (n)| +\sum_{l=1}^j B_{M_n +l}}.
$$
}
\sid{By induction over $n$,} we have
\[
\sum_{v\in \bfV(n)}I_v(n) = 1+M_n = \sum_{v\in \bfV(n)}O_v(n), 
\]
which agrees with \eqref{eq:sum_POPA}.

\sid{For fixed $n\ge 0$, $(\bfI(n), \bfO(n))$ represent in- and
  out-degrees of nodes in $\bfG(n)$ after creation of
  $|\bfE(n)|=1+M_n$ edges and
  $1+\sum_{i=1}^{M_n} B_i$ nodes}
and
\sid{for $k\in \{M_n+1,\ldots M_{n+1}\}$, the choice variables
  $\mathcal{J}_k$, $\widetilde{\mathcal{J}}_k$ have unchanging
  probabilities since $(\bfI(n),\bfO(n))$ only get updated after an
  additional $\Delta M_{n+1}$ edges are added to $\bfG(n)$.  }

Following a similar argument as in Section~\ref{subsec:PA_setup}, we see that
for $n\ge 0$, the transition probability from 
$\bigl(\bfI(n), \bfO(n)\bigr)$ to $\bigl(\bfI(n+1), \bfO(n+1)\bigr)$ 
agrees with the attachment probabilities in a PA model with Poisson measurement given in \eqref{eq:POPA_in} and \eqref{eq:POPA_beta}.
Therefore, $\bigl(\bfI(\cdot), \bfO(\cdot)\bigr) $ represents the evolution of in- and out-degrees in a PA model with Poisson measurement.

\subsection{Degree Distribution}\label{subsec:POPA_deg}
In this section, we focus on the in- and out-degree distributions in a PA model with Poisson measurement, and compare them with those in a traditional PA model.

Set also $\mathcal{J}_0 = \widetilde{\mathcal{J}}_0 \equiv 1$, then
we see from the construction in Section~\ref{subsec:def_POPA} that
\begin{align*}
&\bigl(\bfI(n), \bfO(n)\bigr) = 
\left(\left(\sum_{k=1}^{M_n}\ind_{\{\mathcal{J}_k = v\}}, \sum_{k=1}^{M_n}\ind_{\{\widetilde{\mathcal{J}}_k = v\}}\right)_{v\in\bfV(n-1)},\underbrace{(0,1),\ldots, (0,1),}_{\text{$\sum_{k=1}^{\Delta M_{n}}B_{M_n+k}$ entries}}(0,0), \ldots\right)\\
&=\left(\left\{\left(\sum_{k=M_i}^{M_n}\ind_{\{\mathcal{J}_k = v\}}, \sum_{k=M_i}^{M_n}\ind_{\{\widetilde{\mathcal{J}}_k = v\}}\right)_{v\in\bfV(i)\setminus\bfV(i-1)}\right\}_{0\le i\le n-1},\underbrace{(0,1),\ldots, (0,1),}_{\text{$\sum_{k=1}^{\Delta M_{n}}B_{M_n+k}$ entries}} (0,0), \ldots\right).
\end{align*}
Similar to the induction argument in \eqref{eq:Js_indep}, we have that given $\{M_k: 1\le k\le n\}$
and $\{B_k: 1\le k\le M_n\}$, $\{\mathcal{J}_k: 1\le k\le M_n\}$ and 
$\{\widetilde{\mathcal{J}}_k: 1\le k\le M_n\}$ are conditionally independent.

With a martingale argument, we have the following proposition which summarizes the asymptotic behavior of the in- and out-degrees for a fixed node $v\in \bfV(n)$,
with proof collected in Section~\ref{subsec:pf2}.
\bpr\label{prop:conv_deg_PGPA}
Suppose $\bigl(\bfI(n), \bfO(n)\bigr)$ are as defined in Section~\ref{subsec:def_POPA}.
Then for $v\in \bfV(i)\setminus \bfV(i-1)$, $i\ge 0$, there exists random variables, $\xi_i$, and $\widetilde\xi_i$, such that as $n\to\infty$,
\begin{align*}
\frac{I_v(n)}{\prod_{k=0}^{n-1}\left(1+\frac{\lambda+1}{M_k+1+\deltain|\bfV(k)|}\right)}
&\convas \xi_i,\qquad
\frac{O_v(n)}{\prod_{k=0}^{n-1}\left(1+\frac{(\lambda+1)(1-p)}{M_k+1+\deltaout|\bfV(k)|}\right)}
\convas \widetilde\xi_i.
\end{align*}
Further, 
using the notation $X_n\asymp n^a$, $a>0$, to denote the scenario where
\[
\limsup_{n\to\infty} X_n/n^a <\infty,\quad \text{and}\quad \limsup_{n\to\infty} n^a/X_n <\infty,
\]
 we have
\begin{align*}
\prod_{k=0}^{n-1}&\left(1+\frac{\lambda+1}{M_k+1+\deltain|\bfV(k)|}\right)\asymp n^{1/(1+\deltain p)},\\
\prod_{k=0}^{n-1}&\left(1+\frac{(\lambda+1)(1-p)}{M_k+1+\deltaout|\bfV(k)|}\right)\asymp n^{(1-p)/(1+\deltaout p)}.
\end{align*}
\epr

Here $\left\{\left(\mathcal{J}_{k} , \bigwedge_{v\in \bfV(n)}\frac{\tau_{k,v}}{I_v(n)+\deltain}\right) : k\ge 1\right\}$ can be embedded into a similar birth-immigration framework as in Section~\ref{subsec:PA_dist} (so are $\left\{\left(\widetilde{\mathcal{J}}_{k} , \bigwedge_{v\in \bfV(n)}\frac{\widetilde{\tau}_{k,v}}{O_v(n)+\deltain}\right) : k\ge 1\right\}$), but the difference in the embedding procedure is that 
after observing the values of the $|\bfV(n)| = 1+\sum_{k=1}^{M_n} B_k$ birth-immigration processes,
their transition rates must remain unchanged
 during the occurrence of next $\Delta M_{n+1}$ jumps among the $|\bfV(n)|$ processes. 
We refer to such continuous time period under the embedding framework as the \emph{observation period} (see the proof of Theorem~\ref{thm:IO_Pois} in Section~\ref{subsec:pf2} for details).

\begin{figure}[h]
\centering
\includegraphics[scale=.5]{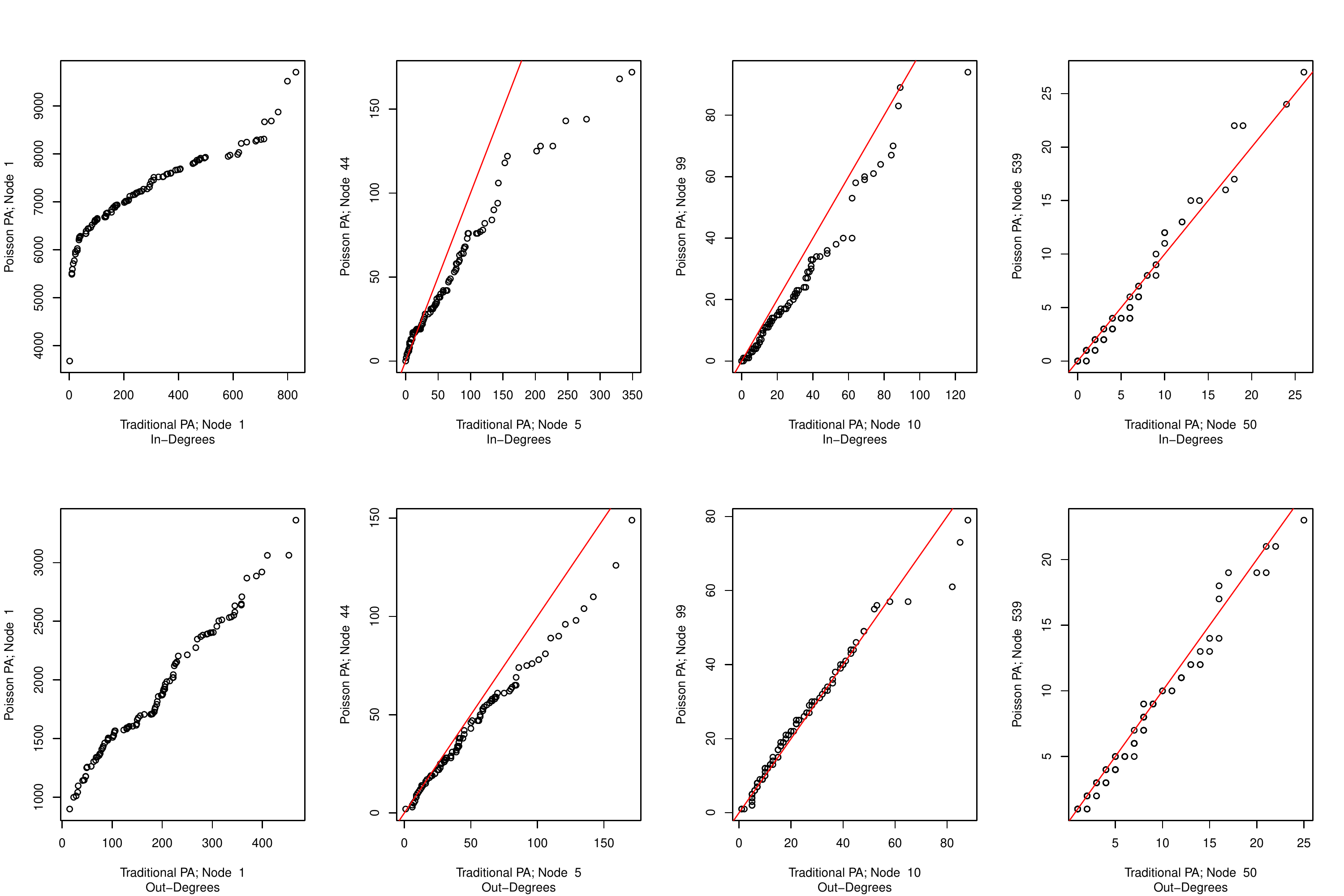}
\caption{QQ plots that compare the in- and out-degree distribution between Node $i$ in a traditional PA model and Node $(\lambda+1)(i-1)+1$ in a Poisson PA model, with $i = 1,5,10,50$. The red reference line is $45^{\circ}$-line. 
When comparing the in- and out-degree distributions for Node 1 in both models, the red line is not displayed since the discrepancy is too large.
We see that overall,
the difference in the degree distribution is large for nodes created earlier in the network, but it vanishes for nodes
added later in the network.}\label{fig:deg_compare}
\end{figure}

Heuristically, for nodes created earlier in the network, the
discrepancy in the degree distribution between the traditional and
Poisson PA models is large. 
For example, the first node in the Poisson PA model has a large
in-degree, as all of the first $M_1$ edges have to point to Node 1 and
there \sid{are} no other existing nodes in the model, which makes Node 1 more
advantageous compared to nodes added later. The first node in a
traditional PA model, however, does not have such behavior.  
For nodes created later in the Poisson PA model, 
the difference between the two models becomes negligible (cf. Figure~\ref{fig:deg_compare}). 

We now give the theoretical justification for the small difference between
 the degree distribution in a PA model with Poisson measurement with that in a traditional PA model for nodes added
 later into the network.
\bthe\label{thm:IO_Pois}
Consider $\bigl(\bfI(n), \bfO(n)\bigr)$ as defined in Section~\ref{subsec:def_POPA}, and we have
for $v\in \bfV(i)\setminus\bfV(i-1)$, $n\epsilon \le i\le n$, $\epsilon>0$, 
\begin{align*}
\PP\left(I_v(n) = m, O_v(n) = l\right) &= \PP\left(\barI_{|\barbV(i)|}(2n) = m, \barO_{|\barbV(i)|}(2n) = l\right)+o\left(i^{-1}\right),
\end{align*}
which further implies as $n\to\infty$,
\beqq\label{eq:approx_pij}
\frac{1}{n}\sum_{i=1}^n  \PP\left(I_{|\bfV(i)|}(n) = m, O_{|\bfV(i)|}(n) = l\right)
\to  p_{m,l}.
\eeqq
In fact, following a similar embedding argument as in Theorem~\ref{thm:PA_BI},
we have, for the PA model with Poisson measurement,
\beqq\label{eq:POPA_pij}
\frac{1}{|\bfV(n)|}\sum_{v\in\bfV(n)}  \ind_{\left\{I_{v}(n) = m, O_{v}(n) = l\right\}}
\convp  p_{m,l},\qquad (m,l)\in \mathbb{N}\times\mathbb{N}_{>0},
\eeqq
where the explicit form of $p_{m,l}$ is the same as in \eqref{eq:PA_pij}.
\ethe

In Figure~\ref{fig:deg_compare}, we provide a numerical comparison between specific nodes in the two PA models considered in this paper. For both models, we choose $(p, \deltain,\deltaout) = (0.2,1,1)$, and generate 100 replications for each model.
When simulating the Poisson PA model, we set $\lambda = 10$ and
the number of unit-shifted Poisson random variables to be
 $n = 2000$. Let the number of edges equal to $2000$ when simulating a traditional PA model. 
Since for each $1\le i\le n$,
the in- and out-degree distributions of $I_v(n)$ are identical for all $v\in \bfV(i)\setminus \bfV(i-1)$,
 we compare the in- and out-degree distributions of Node $i$ in a traditional PA model with Node $(\lambda+1)(i-1)+1$ in a Poisson PA model, where we pick $i = 1,5,10,50$. The comparison is done through QQ-plots and the red reference line is $45^{\circ}$-line.

From Figure~\ref{fig:deg_compare}, we see a huge discrepancy in both in- and out-degree distributions between two models, especially for Node 1 where the $45^{\circ}$ reference line is not even displayed. Such difference then tapers off when it comes to nodes added later into the network, and we observe points from a QQ-plot line up closely with the $45^{\circ}$ reference line.


\section{Model fitting}\label{sec:fit}
\subsection{Background}
We begin this section with some useful results and estimation methods that will be used to fit the directed
PA model with Poisson measurement to real data.

By the formula for $p_{m,l}$ in \eqref{eq:PA_pij}, we have the following corollary which gives 
the marginal power-law behavior for the tail distribution of both in- and out-degrees. 
\begin{Corollary}
Consider the in- and out-degree sequence $\bigl(\bfI(n), \bfO(n)\bigr)_v$ in a directed PA model with Poisson measurement, then
for $m\ge 0$, $l\ge 1$,
\begin{align*}
\frac{1}{|\bfV(n)|}\sum_{v\in \bfV(n)}\ind_{\{I_v(n) = m\}}&\convp \int_0^1\PP\left(Z_{\deltain}\left(t^{1/(1+\deltain p)}\right)=m\right)\mathrm{d}t=:\pin_m,\\
\frac{1}{|\bfV(n)|}\sum_{v\in \bfV(n)}\ind_{\{O_v(n) = l\}}&\convp \int_0^1\PP\left(1+\widetilde{Z}_{1+\deltaout}\left(t^{(1-p)/(1+\deltaout p)}\right)=l\right)\mathrm{d}t=:\pout_l,
\end{align*}
where 
\begin{align}
\pin_m 
\sim (1+\deltain p)\frac{\Gamma(1+\deltain(1+p))}{\Gamma(\deltain)} m^{-(1+\ain)},&\qquad m\to\infty,\nonumber \\
\pout_l \sim \left(\frac{1+\deltaout p}{1-p}\right) \frac{\Gamma\left(1+\frac{1+\deltaout}{1-p}\right)}{\Gamma(1+\deltaout)}
l^{-(1+\aout)},&\qquad l\to\infty,\nonumber\\
\ain = {1+\deltain p},\qquad \aout =\frac{1+\deltaout p}{1-p} .&\label{eq:iota}
\end{align}
\end{Corollary}

One common way in the extreme value theory
to estimate tail indices $(\ain, \aout)$ is to use the Hill estimator
\cite{hill:1975,resnickbook:2007, dehaan:ferreira:2006}.
For non-iid network data, the use of Hill estimator requires justification.
With the distributional results in Theorem~\ref{thm:IO_Pois} available, we presume that the proof machinery in 
\cite{wang:resnick:2018, wang:resnick:2019} to obtain the consistency of Hill estimator is applicable to
data generated from the directed PA model with Poisson measurement.
Hence, we proceed to estimate $\ain$ and $\aout$ by the corresponding
Hill estimator. 

Here we give the estimator for $\ain$ and that for $\aout$ follows in an analogous way.
 Let $I_{(1)}(n) \ge \ldots \ge I_{(|\bfV(n)|)}(n)$ be the decreasing order
statistics of $I_v(n)$, $v\in \bfV(n)$. The Hill estimator
$\hatain(k_n) $ based on $k_n$ largest degrees is
\begin{equation}\label{e:Hill}
	\hatain (k_n) 
	= \left(\frac{1}{k_n} \sum_{j=1}^{k_n} \log\frac{I_{(j)}(n)}{I_{(k_n+1)}(n)}\right)^{-1},
\end{equation}
where $\{k_n\}$ is an intermediate sequence satisfying $k_n \to
\infty$ and $k_n/n \to 0$, as $n\to\infty$. 

To select $k_n$ in practice, the authors of \cite{clauset:shalizi:newman:2009} have proposed computing
 the KS distance between the
empirical distribution tail of the upper $k$ observations and the power-law
distribution with index $\hatain(k)$:
\[
\mathcal{D}_{k}:=\supy \left|\frac{1}{k} \sum_{j=1}^{k}{\bf1}_{\{I_{(j)}(n)/I_{(k+1)}(n)>y\}}-y^{-\hatain(k)}\right|, \quad 1\le k\le n-1.
\]
Then the optimal $k^*_n$ is  the one that minimizes the KS distance:
$$
k^*_n := \argmin_{1\le k\le n} \mathcal{D}_{k},
$$
and the tail index is estimated by $\hatain(k^*_n)$. We refer to the above procedure as the {\it minimum distance method}.
It has been widely used by
data repositories of large network datasets such as KONECT (\url{http://konect.cc/}) \cite{kunegis:2013} and is realized in the R-package {\it poweRlaw\/} \cite{gillespie:2015} as well as the \verb6plfit6 function available at \url{http://tuvalu.santafe.edu/~aaronc/powerlaws/plfit.r}.
For problems with this method, we direct interested readers to \cite{drees:janssen:resnick:wang:2020}.

\medskip

When assessing the goodness of fit, in addition to the comparison among the empirical tail distributions of the in- and out-degrees, another
important way is to inspect the angular density plot, which measures the dependence between in- and out-degrees in a network.
The limit angular density for a directed PA model with Poisson measurement is specified as below.
\bpr\label{prop:ang}
For a pair of random variables $(\mathcal{I},\mathcal{O}) \in \mathbb{N}\times \mathbb{N}_{>0}$ with a joint pmf as given in the right hand side of
\eqref{eq:PA_pij}, consider the following transformation:
\[
(\mathcal{I},\mathcal{O})\mapsto \left(\frac{\mathcal{I}^a}{\mathcal{I}^a +\mathcal{O} }, \mathcal{I}^a + \mathcal{O}\right) =: (\theta, R),\qquad a :=\frac{\iota_\text{in}}{\iota_\text{out}}.
\]
Then using \eqref{eq:PA_pij}, the conditional distribution of $\theta$ given $R>r$ as $r\to\infty$ converges to a distribution function $F$
on $[0,1]$ with density
\begin{equation}\label{eq:ang_PA}
f(\theta) \propto \frac{p}{\deltaout} \theta^{{\deltain}/{a}-1}(1-\theta)^{\deltaout}\int_0^\infty t^{a-1+\iota_\text{in}+\deltain + a\deltaout}
e^{-t\theta^{1/a}-t^a (1-\theta)}\mathrm{d}t,\quad \theta\in[0,1].
\end{equation}
\epr
\begin{proof}
Since the limit density given in Corollary~\ref{cor:PA_count} is the same as the traditional PA model, 
the proof is analogous to that in \cite[Corollary 4.1 and Section 4.1.2]{wan:wang:davis:resnick:2017}, and we omit it here.
\end{proof}

We refer to \eqref{eq:ang_PA} as the limit angular density that measures the asymptotic dependence structure between the in- and out-degrees
in a directed PA network model with Poisson measurement.
To plot the estimated angular density (for example, as in Figure~\ref{fig:fb_PoisPA} (right panel)), 
we first approximate $a$  by $\widehat{a} = \widehat{\iota}_\text{in}/\widehat{\iota}_\text{out}$.
Then the distribution $F$ is estimated via the distribution
of the sample angles 
$$\theta_v(n):= \frac{\left(I_v(n)\right)^{\widehat{a}}}{\left(I_v(n)\right)^{\widehat{a}}+ O_v(n)}, \qquad v\in \bfV(n),$$
for which $R_v(n):=\left(I_v(n)\right)^{\widehat{a}}+ O_v(n)$ exceeds a
large threshold $r$ (chosen to be the $99.5\%$-percentile of $\{R_v(n): v\in \bfV(n)\}$ for all cases considered in this paper).  This is the
POT (Peaks Over Threshold) methodology {commonly employed} in extreme
value theory \cite{coles:2001}.  
This estimation procedure is similar to the extreme value estimation method proposed in \cite{wan:wang:davis:resnick:2017b}, 
where a polar coordinate
transformation with the $L_2$-norm is used to derive the angular density.

Recall the observation from Figure~\ref{fig:deg_compare} that in a PA model with Poisson measurement,
the first node may
have extremely large in- and out-degrees, especially when the Poisson rate parameter $\lambda$ is large, as
it is advantageous in attracting new edges, thus creating a situation where model evolution is slow to forget initial conditions.
To overcome this issue, we scale the estimated $\lambda$ to a smaller time unit.
When the timestamp information is coarse and we are not able to obtain the hourly estimated 
$\lambda$ directly, 
the Poisson assumption allows us to scale the daily estimated $\lambda$ by 24 to get the hourly estimate.

\subsection{Facebook wall posts}
Recall the three plots in Figure~\ref{fig:pois}, and
we here only consider the edges (i.e. wall posts) created from 2007-04-08 to 2008-05-31 ($380,014$ edges in total) , during which 
all three plots remain relatively stable, and discard
the network evolution prior to 2007-04-08. 

\begin{figure}[h]
\centering
\includegraphics[scale=.35]{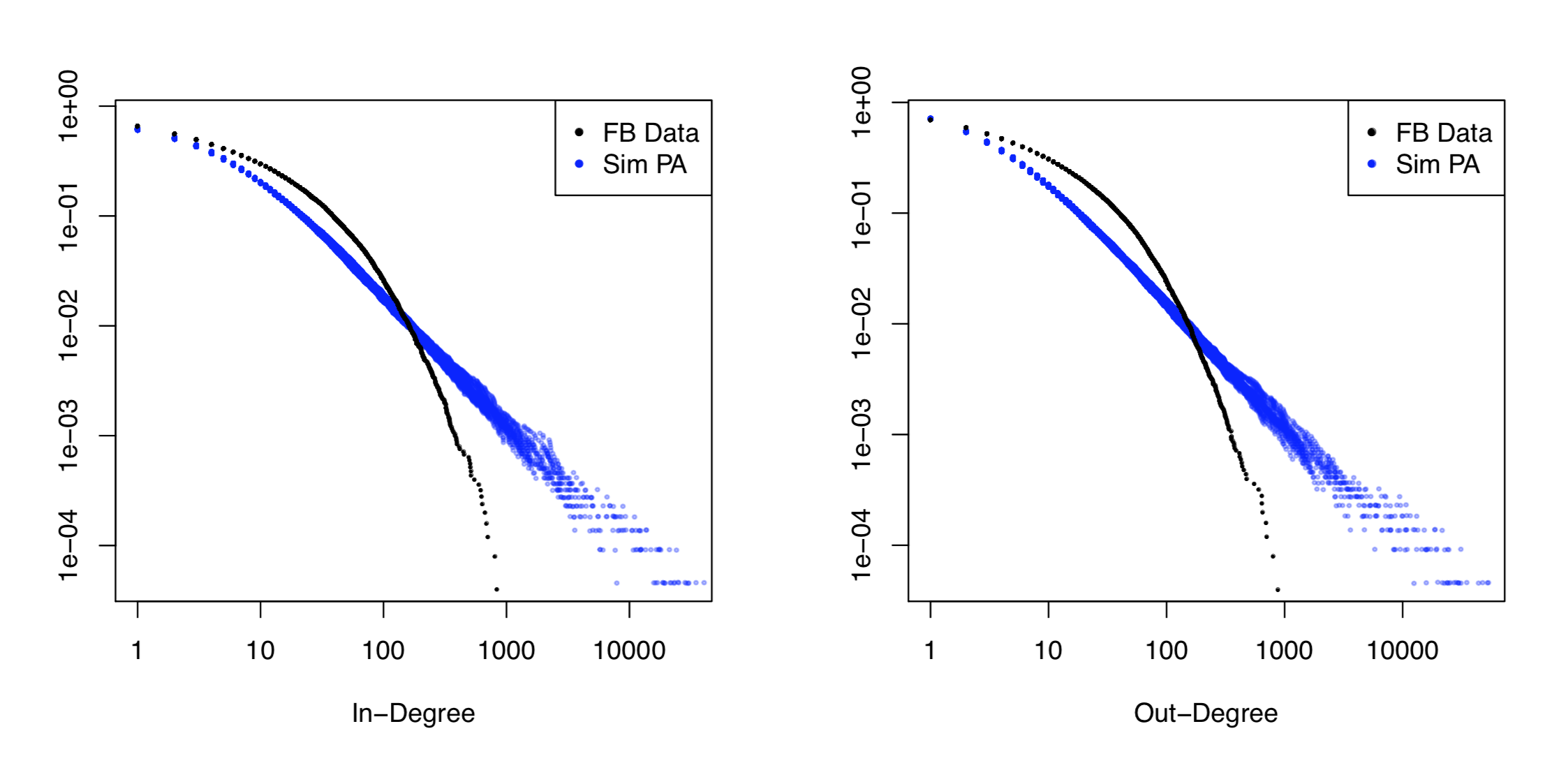}
\caption{A traditional PA model with two scenarios is fitted to the Facebook data using the estimation method in \cite[Section 5]{wan:wang:davis:resnick:2017}. 
Empirical tail distribution of the in- and out-degrees from: (1) Facebook wall post data from 2007-04-08 to 2008-05-31 (black dots); 
(2) 20 replications of the simulated PA data (blue dots) using $\hat{\boldsymbol{\theta}}_\text{MLE}^{\,\text{Fb}}$. Here the MLE method does not provide a good fit to the degree distribution.}\label{fig:fb_mle}
\end{figure}

We first fit a two-scenario traditional PA model (i.e. the one constructed in Section~\ref{subsec:PA_setup}) to the Facebook data using the MLE method proposed in \cite{wan:wang:davis:resnick:2017}.
Applying the MLE estimation algorithm (cf. \cite{wan:wang:davis:resnick:2017}) gives 
\[
\hat{\boldsymbol{\theta}}_\text{MLE}^{\,\text{Fb}} := (\hat{p},
\hat\delta_\text{in}, \hat\delta_\text{out}) 
= (0.057, 3.992, 1.867).
\]
Using $\hat{\boldsymbol{\theta}}_\text{MLE}^{\,\text{Fb}}$, 
we generate 20 independent replications of the traditional PA model using the simulation algorithm in \cite{wan:wang:davis:resnick:2017} 
(each replication contains $380,014$ edges),
and the empirical tail distributions of the in- and out-degrees from the 20 replications are plotted using blue dots in Figure~\ref{fig:fb_mle}. 
In comparison, the empirical tail in- and out-degree distributions from the Facebook data are marked as black dots in Figure~\ref{fig:fb_mle}. 
Due to the huge discrepancy in both tail distributions, we see that fitting a traditional PA model using the MLE method does not provide a good fit for the Facebook data. 

We now fit the newly proposed PA model with Poisson measurement to the Facebook data.
We estimate $p$ by
\[
\widehat{p}^{\,\text{Fb}}= \frac{\text{Total Number of Nodes}}{\text{Total Number of Edges}} = \frac{25085}{380014} \approx 0.066.
\]
Then we use the minimum distance method proposed to obtain estimates for $(\ain, \aout)$:
\beqq
\widehat{\iota}^\text{Fb}_\text{in} = 2.41,\qquad \widehat{\iota}^\text{Fb}_\text{out} = 2.67.
\eeqq\label{eq:iota_hat}
Combining \eqref{eq:iota_hat} with \eqref{eq:iota} and $\widehat{p}^{\,\text{Fb}}$ gives
\[
\widehat{\delta}^{\,\text{Fb}}_\text{in} = 21.42,\qquad \widehat{\delta}^{\,\text{Fb}}_\text{out} = 22.66.
\]
During the period from 2007-04-08 to 2008-05-31 ($n=420$ days), we 
estimate the daily Poisson rates by taking the reciprocal of the averaged inter-events times within a day 
(with timestamps generated during 1-8 AM excluded), and average all 420 daily estimates to obtain
$\widehat\lambda_d^{\,\text{Fb}} = 791.17$, the Poisson rate parameter in the PA model with Poisson measurement.

With $\widehat\lambda_d^{\,\text{Fb}} = 791.17$, when we simulate the Poisson PA model, there are approximately $800$ edges
linking to the Node 1 at the first step. In this sense, the first few nodes created at the beginning of the network will 
distort the degree distribution. Rescaling the estimated daily Poisson rate to the hourly rate and elongating time scale $n$ to be 
the total number of hours over which the network has evolved provide remedies to such problems.
Note also that $\widehat\lambda_d^{\,\text{Fb}}$ is calculated after expected sleeping hours (1--8 AM) are excluded. 
Therefore, after rescaling, the
 estimated hourly Poisson rate becomes
 \[
 \widehat{\lambda}_h^{\,\text{Fb}} = \widehat\lambda_d^{\,\text{Fb}}/17 = 46.54,
 \]
and $n = 420\times 17 = 7140$.

\begin{figure}[h]
\centering
\includegraphics[scale=.35]{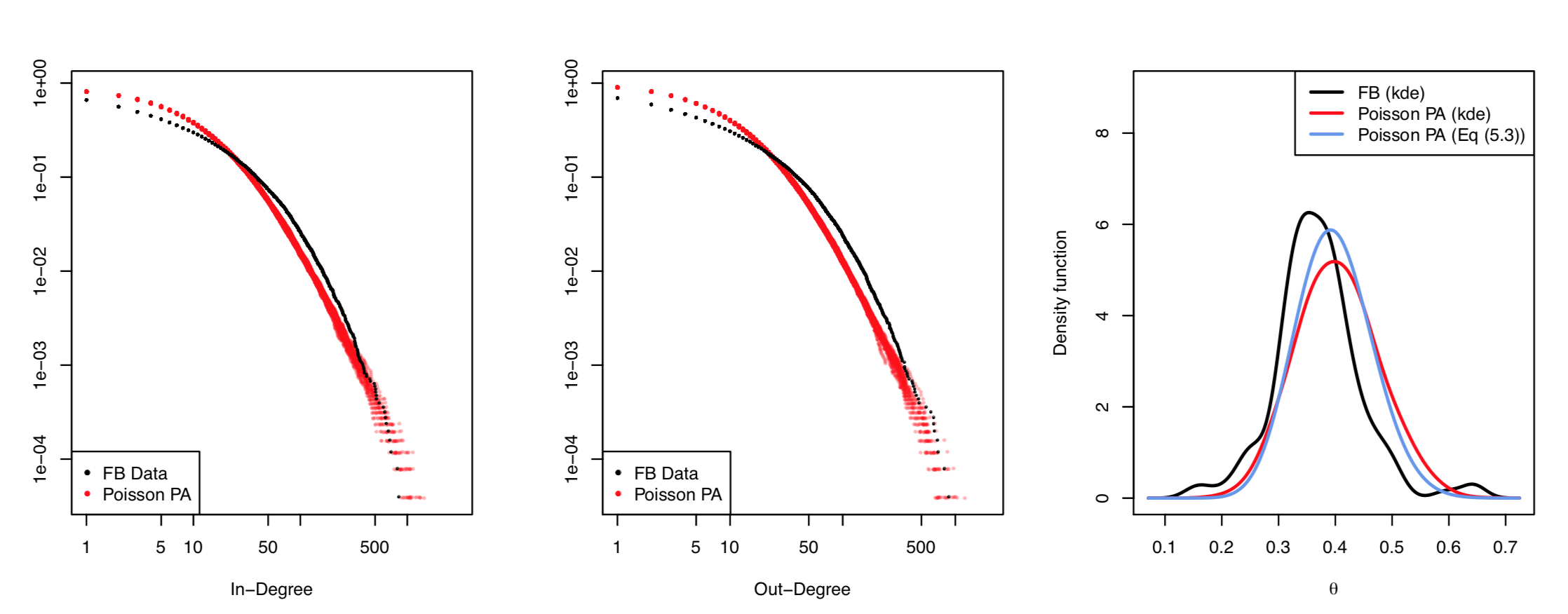}
\caption{With $(\widehat{\lambda}_h^{\,\text{Fb}},\widehat{p}^{\,\text{Fb}}, \widehat{\delta}_\text{in}^{\,\text{Fb}}, \widehat{\delta}_\text{out}^{\,\text{Fb}})
 = (46.54, 0.066, 21.42, 22.66)$
and 
 $n^{\,\text{Fb}} = 7140$, 
 20 independent replications of 
the directed PA model with Poisson measurement are simulated. Left: Empirical tail distribution of in-degrees. 
Middle: Empirical tail distribution of out-degrees. Right: Estimated angular densities using the non-parametric method (red and black) as well as \eqref{eq:ang_PA} and parameter estimates (blue). Red labels correspond to the simulated data, while
the black ones represent the Facebook data. Compared with Figure~\ref{fig:fb_mle}, the PA model with Poisson measurement
provides a better fit.}\label{fig:fb_PoisPA}
\end{figure}

With $\widehat{\boldsymbol{\theta}}_{\text{Fb}} := (\widehat{\lambda}_h^{\,\text{Fb}},\widehat{p}^{\,\text{Fb}}, \widehat{\delta}_\text{in}^{\,\text{Fb}}, \widehat{\delta}_\text{out}^{\,\text{Fb}})$ available and $n^{\,\text{Fb}} = 7140$, we simulate 20 independent replications of 
the directed PA model with Poisson measurement described in Section~\ref{subsec:POPA}.
The empirical tail distributions of the in- and out-degrees from the 20 replications are plotted in the left and middle panels
in Figure~\ref{fig:fb_PoisPA} using red dots. Compared with the degree distributions in Figure~\ref{fig:fb_mle}, 
the PA model with Poisson measurement apparently provides a better fit, though slight discrepancy still exists. 

The right panel in Figure~\ref{fig:fb_PoisPA} compares the estimated angular density, using the results in Proposition~\ref{prop:ang}.
With $\{(I_v(n), O_v(n): v\in \bfV(n)\}$ observed from the Facebook data, 
we calculate $\theta_v(n)$ and use the \verb6kde6 function in the \verb6R6 package \verb6ks6
to get the estimated density of $\{\theta_v(n)$: $v\in V(n)\}$, thus generating the black curve in the right panel of Figure~\ref{fig:fb_PoisPA}.
For the 20 simulated Poisson PA networks, the estimated angular densities are calculated in the same way, but due to the variation across different replications,
we only report the averaged angular density estimates in the right panel of Figure~\ref{fig:fb_PoisPA} (the red curve).

Figure~\ref{fig:fb_PoisPA} reveals that the in- and out-degrees in the Facebook data are asymptotically dependent, and 
the mode is around $\theta = 0.35$.
The angular density based on the simulated data from the PA model with Poisson measurement
is unimodal with a mode around 0.4.
Further adjusting the estimated Poisson rate with a narrower time window (e.g. from hourly to 30-min)
pushes the mode of the estimated angular density closer to 0.35.
We plot the estimated asymptotic angular density given in \eqref{eq:ang_PA} (with
$\widehat{\boldsymbol{\theta}}_{\text{Fb}}$ plugged in) as the blue curve in the right panel of Figure~\ref{fig:fb_PoisPA}.

One important message here is applying the MLE method to the Facebook data gives a much poorer fit than fitting our modified PA model
with Poisson measurement. One possible explanation is that the MLE method is less robust 
to data corruption and model mis-specification, compared to the estimation approach using extreme value theory
to first estimate $\ain$ and $\aout$ (see \cite{wan:wang:davis:resnick:2017b} for 
more examples).

\subsection{Slashdot}
In this section, we discuss another social network dataset, \emph{Slashdot}, which is mentioned in Section \ref{sec:intro_data}. 
One special feature of the Slashdot data is the coarse timestamp information, where edge creation times are recorded to the nearest minute.
So several new edges can be added with the same timestamp. Meanwhile, in this dataset, it is not necessarily true that
a node labeled with a smaller number is created at an earlier time point. Hence, fitting a traditional PA model using the MLE method
is not applicable. 

\begin{figure}[h]
\centering
\includegraphics[scale=.6]{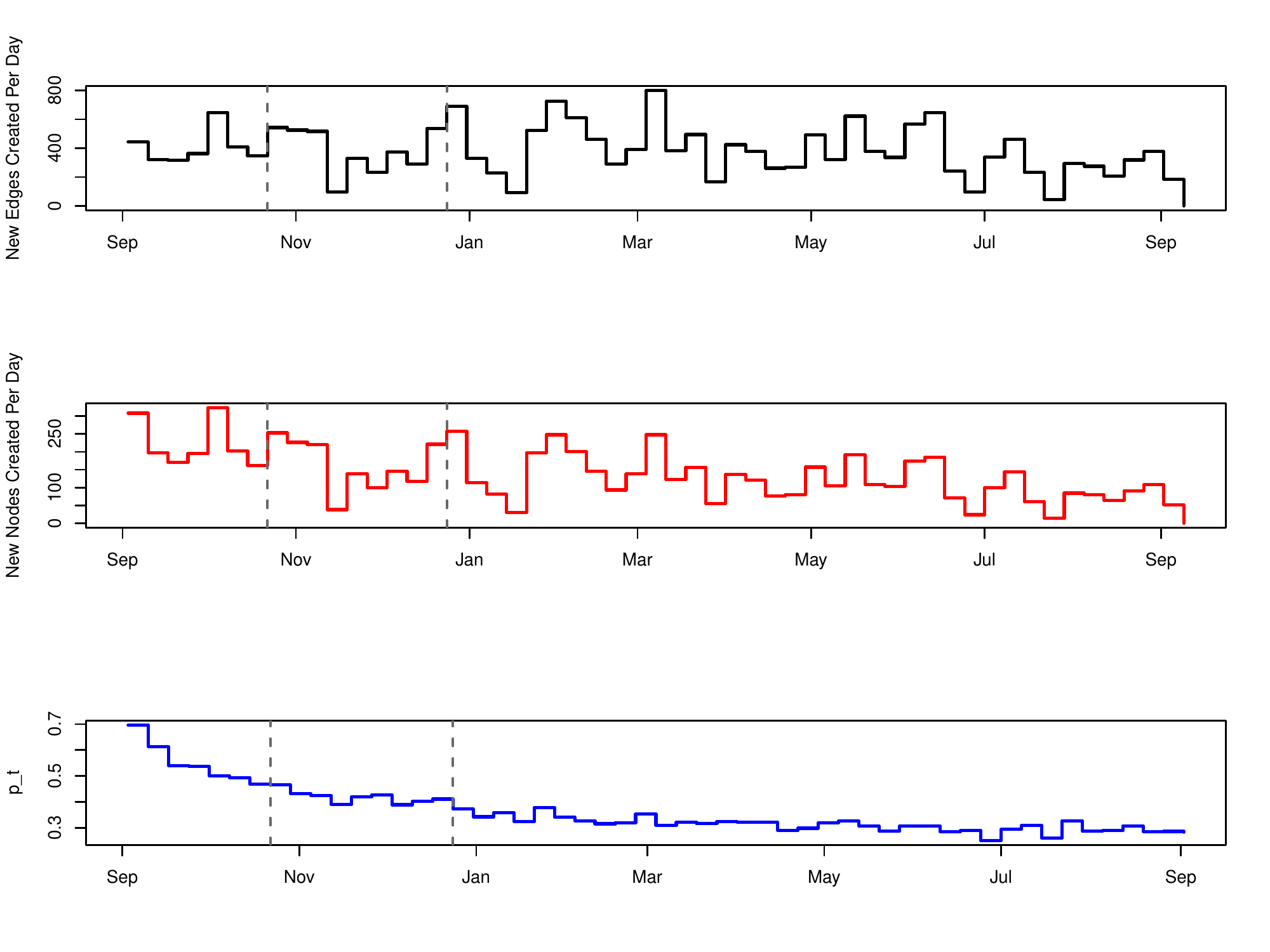}
\caption{Daily new edge counts (top) and node counts (middle) (averaged by week) from 2005-09-01 to 2006-09-02. The
bottom panel gives the ratio of new node counts over new edge counts. 
Within the third time interval (from 2005-12-18 to 2006-09-02), all three plots remain relatively stable.}\label{fig:sd_pois}
\end{figure}

Learning from the findings in \cite{wang:resnick:2019b}, we assume that the edge creation process follows an NHPP
which has a constant rate within a day but varying rates from day to day.
With only coarse timestamp information, we estimate the Poisson rates of edge and node creation processes by
averaging new edges and nodes created per day over non-overlapping weekly intervals, respectively. The estimation results are plotted in 
Figure~\ref{fig:sd_pois}.

The top and middle panels in Figure~\ref{fig:sd_pois} report the 
averaged daily new edge and node counts over non-overlapping weekly intervals from 2005-09-01 to 2006-09-02, respectively.
The bottom panel gives the ratio of new node counts over new edge counts, $\widehat{p}_t$.
The bottom panel shows a decreasing trend at first and remains relatively stable later. 
We again use the \verb6breakpoints6 function in R's \verb6strucchange6 package to locate the change points of $\widehat{p}_t$.
The two change points detected are marked as gray vertical lines in all panels of Figure~\ref{fig:sd_pois}.
We see that during the last time segment, all three panels display a relatively stable trend.
Hence, we proceed by focusing on the data in the last time segment, namely from
2005-12-18 to 2006-09-02 (259 days).

Again, the two tail indices, $\ain$ and $\aout$, are estimated using the minimum distance method so that
\[
\widehat{\iota}_\text{in}^{\,\text{Sd}} = 2.76,\qquad \widehat{\iota}_\text{out}^{\,\text{Sd}} = 2.06.
\]
We estimate the parameters $(\lambda, p, \deltain, \deltaout)$ from the PA model with Poisson measurement as follows:
\begin{enumerate}
\item $\widehat{\lambda}_d^{\,\text{Sd}} = $ Average daily counts of new edges.
\item $\widehat{p}^{\,\text{Sd}} = \text{Total number of nodes}/ \text{Total number of edges}$.
\item $\widehat{\delta}_\text{in}^{\,\text{Sd}} = (\widehat\iota_\text{in}^{\,\text{Sd}}- 1)/\widehat{p}^{\,\text{Sd}}$.
\item  $\widehat{\delta}_\text{out}^{\,\text{Sd}} = (\widehat\iota_\text{out}^{\,\text{Sd}}\times(1-\widehat{p}^{\,\text{Sd}}) - 1)/\widehat{p}^{\,\text{Sd}}$.
\end{enumerate}
Using the Slashdot data from 2005-12-18 to 2006-09-02, we have
\[
\widehat{\boldsymbol{\theta}}_\text{Sd,d} := (\widehat\lambda_d^{\,\text{Sd}}, \widehat p^{\,\text{Sd}}, \widehat{\delta}_\text{in}^{\,\text{Sd}}, \widehat{\delta}_\text{out}^{\,\text{Sd}})
= (377.29,\, 0.38,\, 4.66, \, 0.73 ).
\]
Similar to the Facebook case, $\widehat\lambda_d^{\,\text{Sd}}$ is so large that simulations will generate at least one node with large in- and out-degrees,
thus distorting the degree distribution. Hence, we rescale $\widehat\lambda_d^{\,\text{Sd}}$ to the hourly Poisson rate:
\[
\widehat{\lambda}_h^{\,\text{Sd}} := \widehat\lambda_d^{\,\text{Sd}}/24 = 15.72,
\]
then
\[
\widehat{\boldsymbol{\theta}}_\text{Sd,h} = (\widehat\lambda_h^{\,\text{Sd}}, \widehat p^{\,\text{Sd}}, \widehat{\delta}_\text{in}^{\,\text{Sd}}, \widehat{\delta}_\text{out}^{\,\text{Sd}})
= (15.72,\, 0.38,\, 4.66, \, 0.73 ).
\]
Note that for the Slashdot data, we do not assume any expected sleeping hours on each day so $\widehat\lambda_d^{\,\text{Sd}}$ is scaled by 24.
With $\widehat{\boldsymbol{\theta}}_\text{Sd,h}$ available and $n^{\,\text{Sd}} = 259\times 24 = 6216$, we simulate 20 independent replications of
the directed PA model with Poisson measurement. 

\begin{figure}
\centering
\includegraphics[scale=.35]{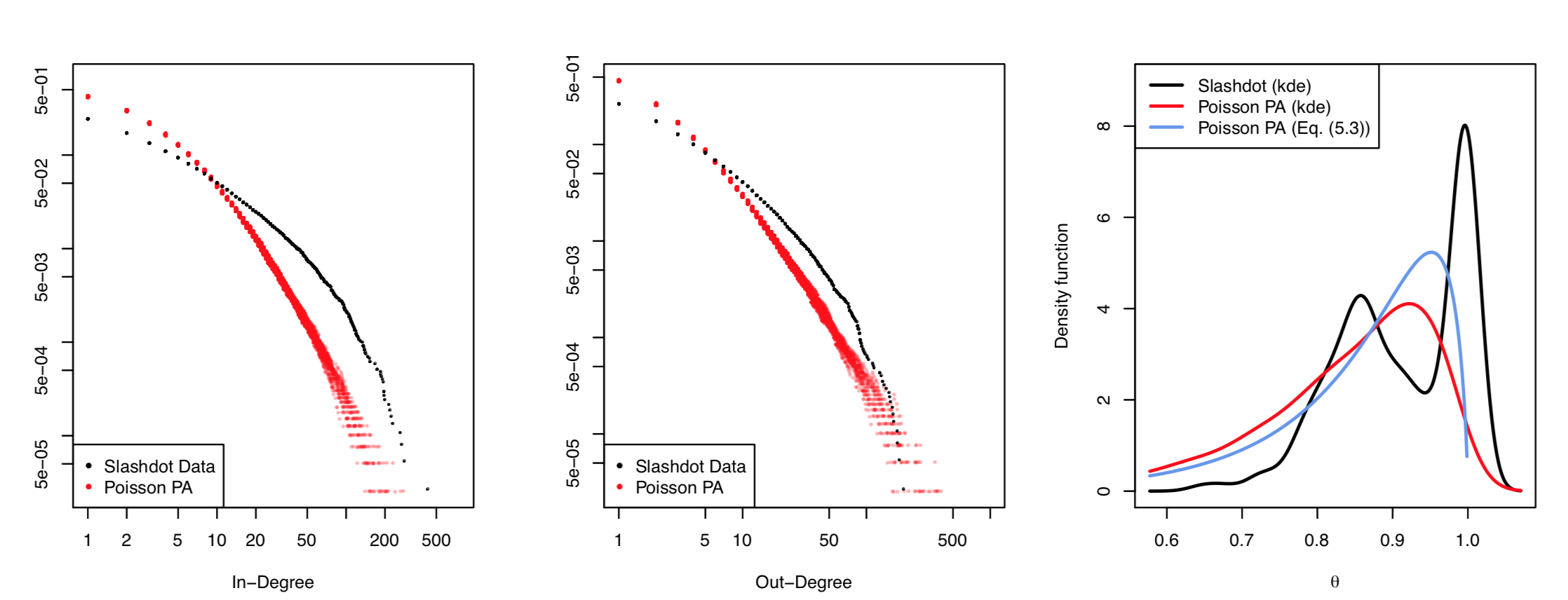}
\caption{With $(\widehat{\lambda}_h^{\,\text{Sd}},\widehat{p}^{\,\text{Sd}}, \widehat{\delta}_\text{in}^{\,\text{Sd}}, \widehat{\delta}_\text{out}^{\,\text{Sd}}, n^{\,\text{Sd}}) 
= (15.72, 0.38,4.66, 0.73,  6216)$, 20 independent replications of 
the directed PA model with Poisson measurement are simulated. Left: Empirical tail distribution of in-degrees. 
Middle: Empirical tail distribution of out-degrees. Right: Estimated angular densities using the non-parametric method (red and black) as well as \eqref{eq:ang_PA} together with parameter estimates (blue). Red labels correspond to the simulated data, while
the black ones represent the Slashdot data. The bimodal shape in the angular density flags discrepancies in dependence modeling.}\label{fig:sd_PoisPA}
\end{figure}

Empirical tail distributions of the in- and out-degrees from the 20 simulated networks (red points) are given in the left and middle panels of Figure~\ref{fig:sd_PoisPA},  and those of the Slashdot data are plotted using black dots. 
The right panel of Figure~\ref{fig:sd_PoisPA} compares the estimated angular densities (using the \verb6kde6 function) based on the in- and out-degrees in Slashdot data (black) 
with the averaged estimated angular densities of the 20 simulated networks (red). 
The blue curve represents the estimated asymptotic angular density \eqref{eq:ang_PA} with $\widehat{\boldsymbol{\theta}}_{\text{Sd},h}$ plugged in.

We observe significant differences in the empirical tail distributions from the left and middle panels of Figure~\ref{fig:sd_PoisPA}.
The estimated angular density in black (based on the Slashdot data) shows a bimodal pattern which the simulated networks fail to catch. 
In fact, if we plot the in- vs out-degrees for all nodes in the Slashdot data (not included here), we see that quite a few nodes have large in-degrees but 0 out-degree. Such unusual pattern explains the higher peak in the right panel of Figure~\ref{fig:sd_PoisPA}, and we speculate these nodes 
are administration accounts that never respond to other users in the network.

\begin{figure}[h]
\centering
\includegraphics[scale=.35]{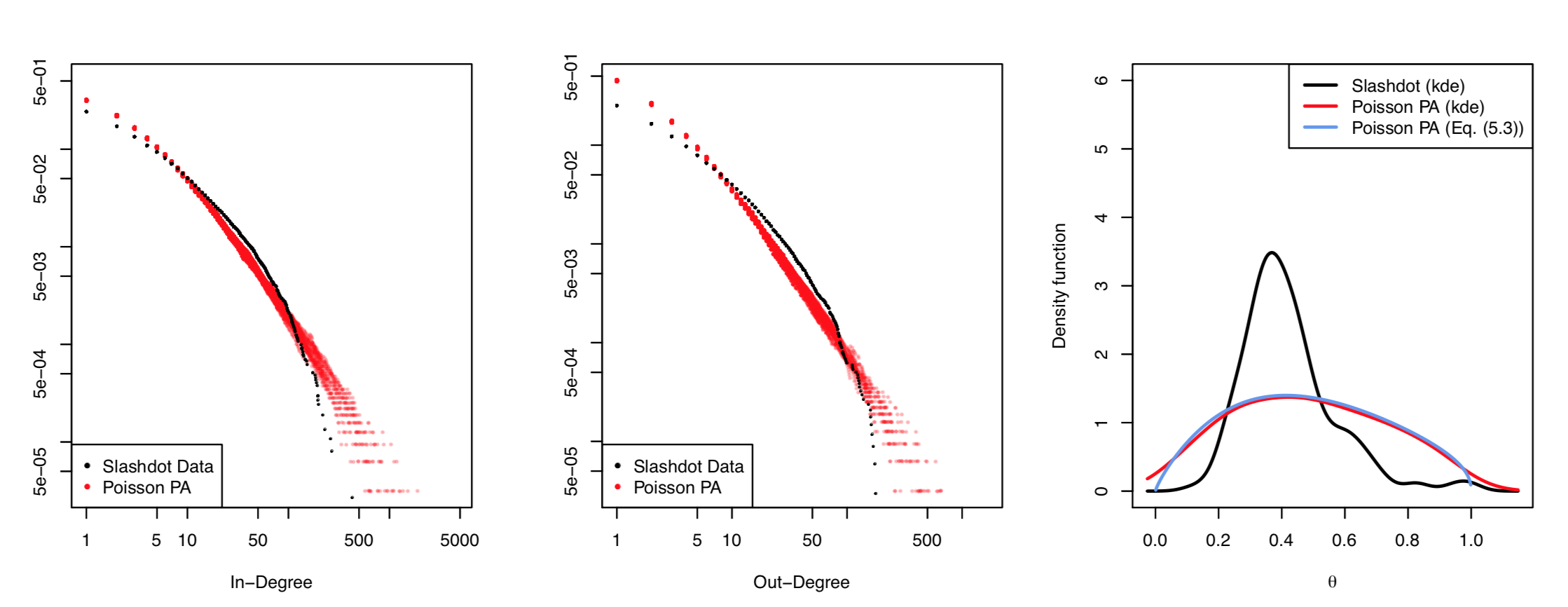}
\caption{With $(\widetilde\lambda_h^{\,\text{Sd}}, \widetilde p^{\,\text{Sd}}, \widetilde{\delta}_\text{in}^{\,\text{Sd}}, \widetilde{\delta}_\text{out}^{\,\text{Sd}}, n^{\,\text{Sd}}) 
= (13.79, 0.34, 1.58, 0.44,  6216)$, 20 independent replications of 
the directed PA model with Poisson measurement are simulated. Left: Empirical tail distribution of in-degrees. 
Middle: Empirical tail distribution of out-degrees. Right: Estimated angular densities using the non-parametric method (red and black) as well as \eqref{eq:ang_PA} together with parameter estimates (blue). Red labels correspond to the simulated data, while
the black ones represent the Slashdot data with administration accounts removed.}\label{fig:sd_PoisPA_new}
\end{figure}

We proceed with the interpretation that all nodes with 0 out-degree and in-degree $\ge 20$ are administration accounts,
After removing all such accounts, 
we refit the directed PA model with Poisson measurement as before. 
The minimum distance method gives $(\widetilde{\iota}_\text{in}^{\,\text{Sd}}, \widetilde{\iota}_\text{out}^{\,\text{Sd}}) = (1.54,1.76)$, then we obtain the following estimates:
\[
\widetilde{\boldsymbol{\theta}}_\text{Sd,h} = (\widetilde\lambda_h^{\,\text{Sd}}, \widetilde p^{\,\text{Sd}}, \widetilde{\delta}_\text{in}^{\,\text{Sd}}, \widetilde{\delta}_\text{out}^{\,\text{Sd}})
= (13.79,\, 0.34,\, 1.58, \, 0.44 ).
\]
Again, we simulate 20 independent PA model with Poisson measurements using $\widetilde{\boldsymbol{\theta}}_\text{Sd,h}$ and 
$n = 6216$, and compare the empirical tail distributions of the in- and out-degrees as well as the angular densities in Figure~\ref{fig:sd_PoisPA_new}.
We observe the discrepancies in the in- and out-degree tail distributions become smaller, and the angular density from the Slashdot data 
becomes unimodal after administration accounts are removed. 

Although the difference between the estimated angular density from the Slashdot data and that from the simulated network still exists,
we note that simply removing nodes with in-degree $\ge 20$ and out-degree equal to $0$ is a rather crude method
to account for the bimodal pattern in the angular density plot given in Figure~\ref{fig:sd_PoisPA} (the right panel).
This may lead to the incorrect asymptotic dependence structure as summarized by the angular density plot. 
We defer studies on the detection of such administration accounts from real data to future work.  

\section{Additional comments}\label{sec:comments}
Note that Chapter 8.2 of \cite{vanderHofstad:2017} considers an undirected PA model where at each step a deterministic $m$ number of edges are added to the network.
In this paper, we intend to consider the theoretical importance of adding a shifted Poisson number of edges. Presumably, the Poisson distribution can be generalized. 
In the absence of statistical evidence of other distributions, we have not
considered this a priority for this project and will consider this in the future.

\section{Proofs}\label{sec:proof}
\subsection{Proofs in Section~\ref{sec:embed_PA}}\label{subsec:pf1}
\begin{proof}[Proof of Theorem~\ref{thm:PA_BI}] 
It remains to show \eqref{eq:PA_BI_dist}. 
According to the embedding framework, we see that
$\left\{(J'_k,\Gamma_{k}-\Gamma_{k-1}) : 1\le k\le n\right\}$
and $\left\{(\widetilde J'_k, \,\widetilde\Gamma_{k}-\widetilde\Gamma_{k-1}): 1\le k\le n\right\}$
are conditionally independent under $\PP^{\{B_k\}_{k=1}^n}$.

By \eqref{eq:Js_indep}, we have for fixed $n$, $m_v\ge 0$, $l_v\ge 1$, $v\in \barbV(n)$,
\begin{align*}
\PP&\left(\bigl(\barI_v(2n),\barO_v(2n)\bigr) = (m_v, l_v), v\in \barbV(n)\right)\\
&= \EE\left[\PP^{\{B_k\}_{k=1}^n}\left(\barI_v(2n)=m_v, v\in \barbV(n)\right)
\PP^{\{B_k\}_{k=1}^n}\left(\barO_v(2n)=l_v, v\in \barbV(n)\right)\right]\\
\intertext{by the embedding results,}
&= \EE\left[\PP^{\{B_k\}_{k=1}^n}\left(\ind_{\{v=1\}} + BI^{(v)}_{\ind_{\{v=1\}}+\deltain}\left(\Gamma_n - \Gamma_{S_v}\right)= m_v, v\in \barbV(n)\right)\right.\\
&\left. \quad\qquad\PP^{\{B_k\}_{k=1}^n}\left(1+\widetilde{BI}^{(v)}_{1+\deltaout}\left(\widetilde\Gamma_n - \widetilde\Gamma_{S_v}\right)=l_v, v\in \barbV(n)\right)\right]\\
&= \PP\left(\left(\ind_{\{v=1\}} + BI^{(v)}_{\ind_{\{v=1\}}+\deltain}\left(\Gamma_n - \Gamma_{S_v}\right),
1+\widetilde{BI}^{(v)}_{1+\deltaout}\left(\widetilde\Gamma_n - \widetilde\Gamma_{S_v}\right)
\right)=(m_v, l_v), v\in \barbV(n)\right).
\end{align*}
\end{proof}
\begin{proof}[Proof of Corollary~\ref{cor:PA_count}]
Consider
\begin{align*}
\frac{1}{|\barbV(n)|}
&\sum_{v\in \barbV(n)} \ind_{\left\{\bigl(\barI_v(2n), \barO_v(2n) \bigr) = (m,l)\right\}}\\ 
&= \frac{1}{|\barbV(n)|}\sum_{v = 2}^{|\barbV(n)|} \ind_{\left\{\bigl(\barI_v(2n), \barO_v(2n) \bigr) = (m,l)\right\}}
+\frac{1}{n}\ind_{\left\{\bigl(\barI_1(2n), \barO_1(2n) \bigr) = (m,l)\right\}},
\end{align*}
and we see that the second term on the right hand side goes to 0 a.s. as $n\to\infty$. It then suffices to consider
\begin{align}\label{eq:PA_counts_new}
\frac{1}{|\barbV(n)|}&\sum_{v = 2}^{|\barbV(n)|} \ind_{\left\{\bigl(\barI_v(2n), \barO_v(2n) \bigr) = (m,l)\right\}},\nonumber\\
\intertext{which for fixed $n$, has the same distribution as (cf. Theorem~\ref{thm:PA_BI}):}
\stackrel{d}{=} 
&\frac{1}{|\barbV(n)|}\sum_{v = 2}^{|\barbV(n)|} \ind_{\left\{\bigl(BI^{(v)}_{\deltain}(\Gamma_n-\Gamma_{S_v}), 1+\widetilde{BI}^{(v)}_{1+\deltaout}(\widetilde{\Gamma}_n - \widetilde{\Gamma}_{S_v}) \bigr) = (m,l)\right\}}\\
= & \frac{1}{|\barbV(n)|}\sum_{k=1}^n B_k\ind_{\left\{\bigl(BI^{(|\barbV(k)|)}_{\deltain}(\Gamma_n-\Gamma_{k}), 1+\widetilde{BI}^{(|\barbV(k)|)}_{1+\deltaout}(\widetilde{\Gamma}_n - \widetilde{\Gamma}_{k}) \bigr) = (m,l)\right\}},\nonumber
\end{align}
\sid{where the coefficient $B_k$ in front of the indicator guarantees we
sum different BI processes inside the indicator.}

Now we divide the quantity in \eqref{eq:PA_counts_new} into different parts:
\begin{align*}
\frac{1}{|\barbV(n)|}&\left[\sum_{v = 2}^{|\barbV(n)|} \ind_{\left\{\bigl(BI^{(v)}_{\deltain}(\Gamma_n-\Gamma_{S_v}), 1+\widetilde{BI}^{(v)}_{1+\deltaout}(\widetilde{\Gamma}_n - \widetilde{\Gamma}_{S_v}) \bigr) = (m,l)\right\}}\right.\\
&\left.\qquad
- \sum_{k=1}^n B_k\ind_{\left\{\left(BI^{(|\barbV(k)|)}_{\deltain}\left(\frac{1}{1+\deltain p}\log\left(\frac{n}{k}\right)\right), 1+ \widetilde{BI}^{(|\barbV(k)|)}_{1+\deltaout}\left(\frac{1-p}{1+\deltaout p}\log\left(\frac{n}{k}\right)\right)\right) = (m,l)\right\}}  \right]\\
&+ \frac{1}{|\barbV(n)|}\sum_{k = 1}^{n}\left[ B_k
\ind_{\left\{\left(BI^{(|\barbV(k)|)}_{\deltain}\left(\frac{1}{1+\deltain p}\log\left(\frac{n}{k}\right)\right), 1+ \widetilde{BI}^{(|\barbV(k)|)}_{1+\deltaout}\left(\frac{1-p}{1+\deltaout p}\log\left(\frac{n}{k}\right)\right)\right) = (m,l)\right\}} \right.\\
& \left. \quad -\, p\PP\left(\left(BI^{(|\barbV(k)|)}_{\deltain}\left(\frac{1}{1+\deltain p}\log\left(\frac{n}{k}\right)\right), 
\right.\right.\right.\\
&\left.\left.\left.\qquad\qquad\qquad\qquad
1+ \widetilde{BI}^{(|\barbV(k)|)}_{1+\deltaout}\left(\frac{1-p}{1+\deltaout p}\log\left(\frac{n}{k}\right)\right)\right) = (m,l)\right)
\right]\\
& + \left[\frac{p}{|\barbV(n)|}\sum_{k = 1}^{n} 
\PP\left(\left(BI^{(|\barbV(k)|)}_{\deltain}\left(\frac{1}{1+\deltain p}\log\left(\frac{n}{k}\right)\right), 
\right.\right.\right.\\
&\left.\left.\left.\qquad\qquad\qquad\qquad\qquad
1+ \widetilde{BI}^{(|\barbV(k)|)}_{1+\deltaout}\left(\frac{1-p}{1+\deltaout p}\log\left(\frac{n}{k}\right)\right)\right) = (m,l)\right)\right.\\
& \left. \quad -\, \int_0^1 \PP\left(\left(Z_{\deltain}\left(t^{\frac{1}{1+\deltain p}}\right), 1+ \widetilde{Z}_{1+\deltaout}\left(t^{\frac{1-p}{1+\deltaout p}}\right)\right) = (m,l)\right)\dd t
\right]\\
& + \int_0^1  \PP\left(\left(Z_{\deltain}\left(t^{\frac{1}{1+\deltain p}}\right), 1+ \widetilde{Z}_{1+\deltaout}\left(t^{\frac{1-p}{1+\deltaout p}}\right)\right) = (m,l)\right)\dd t\\
&=: R_1(n) + R_2(n) + R_3(n) + R_4.
\end{align*}
By \cite[Theorem III.9.1, Page 119]{athreya:ney:2004}, both
\begin{align*}
\Gamma_n &- \sum_{k=0}^{n-1}\frac{1}{k+1+\deltain |\barbV(k)|}\quad
\text{and}\quad
\widetilde\Gamma_n - \sum_{k=0}^{n-1}\frac{1-p}{k+1+\deltaout |\barbV(k)|}
\end{align*}
are $L_2$-bounded martingales with respect to sigma fields
$\mathcal{F}_n'$ \sid{defined in \eqref{e:sigmafield}}
and
$$\sigma\left(\left\{B_k\right\}_{k=1}^n;\left\{\widetilde{BI}^{(v)}_{1+\deltaout}(t-\widetilde{\Gamma}_{S_v}): \widetilde{\Gamma}_{S_v}\le t\le \widetilde{\Gamma}_n\right\}_{v\in\barbV(n)}\right),$$ respectively, which therefore
converge a.s.. Then by \cite[Corollary 2.1(iii)]{athreya:ghosh:sethuraman:2008}, we have
for $\epsilon>0$, 
\begin{align}
\sup_{n\epsilon\le k\le  n}\left|\Gamma_n - \Gamma_{k}
-\frac{1}{1+\deltain p}\log(n/k)\right| &\convas 0,
\label{eq:conv_gamma1}\\
\sup_{n\epsilon\le k\le n}\left|\widetilde\Gamma_n - \widetilde\Gamma_{k}
-\frac{1-p}{1+\deltaout p}\log(n/k)\right|
&\convas 0,
\label{eq:conv_gamma2}
\end{align}
as $n\to\infty$.
Also, note that 
\begin{align*}
|R_1(n)|&\le \frac{1}{|\barbV(n)|}\sum_{k=1}^{n}B_k\left|\ind_{\left\{BI^{(|\barbV(k)|)}_{\deltain}(\Gamma_n-\Gamma_{k})= m\right\}}
\ind_{\left\{1+\widetilde{BI}^{(|\barbV(k)|)}_{1+\deltaout}(\widetilde{\Gamma}_n-\widetilde{\Gamma}_{k}) = l\right\}}\right.\\
 & \left. \qquad \qquad - \ind_{\left\{BI^{(|\barbV(k)|)}_{\deltain}\left(\frac{1}{1+\deltain p}\log(n/k)\right) = m\right\}} 
\ind_{\left\{1+\widetilde{BI}^{(|\barbV(k)|)}_{1+\deltaout}\left(\frac{1-p}{1+\deltaout p}\log(n/k)\right) = l\right\}}  \right|\\
&\le \frac{1}{|\barbV(n)|}\sum_{k=1}^{n}B_k\left|\ind_{\left\{BI^{(|\barbV(k)|)}_{\deltain}(\Gamma_n-\Gamma_{k})= m\right\}} - \ind_{\left\{BI^{(|\barbV(k)|)}_{\deltain}\left(\frac{1}{1+\deltain p}\log(n/k)\right) = m\right\}} \right|\\
&\qquad +\sum_{k=1}^{n}B_k \left|\ind_{\left\{1+\widetilde{BI}^{(|\barbV(k)|)}_{1+\deltaout}(\widetilde{\Gamma}_n-\widetilde{\Gamma}_{k}) = l\right\}} - \ind_{\left\{1+\widetilde{BI}^{(|\barbV(k)|)}_{1+\deltaout}\left(\frac{1-p}{1+\deltaout p}\log(n/k)\right) = l\right\}} \right|.
\end{align*}
Therefore,
\begin{align*}
\EE |R_1(n)| &\le \EE\left(\frac{1}{|\barbV(n)|}\sum_{k=1}^{n}B_k\left|\ind_{\left\{BI^{(|\barbV(k)|)}_{\deltain}(\Gamma_n-\Gamma_{k})= m\right\}} - \ind_{\left\{BI^{(|\barbV(k)|)}_{\deltain}\left(\frac{1}{1+\deltain p}\log(n/k)\right) = m\right\}} \right|\right)\\
 + &\EE\left(\frac{1}{|\barbV(n)|}\sum_{k=1}^{n}B_k\left|\ind_{\left\{1+\widetilde{BI}^{(|\barbV(k)|)}_{1+\deltaout}(\widetilde{\Gamma}_n-\widetilde{\Gamma}_{k}) = l\right\}} - \ind_{\left\{1+\widetilde{BI}^{(|\barbV(k)|)}_{1+\deltaout}\left(\frac{1-p}{1+\deltaout p}\log(n/k)\right) = l\right\}} \right|\right)\\
&=: \EE (R^I_1(n)) + \EE (R^O_1(n)).
\end{align*}
Applying the a.s. convergence results in \eqref{eq:conv_gamma1}, 
\cite[Corollary 3.1]{athreya:ghosh:sethuraman:2008} and the proof machinery for \cite[Theorem 1.2, Page 489--490]{athreya:ghosh:sethuraman:2008}, we obtain that as $n\to\infty$,
\begin{align*}
\lim_{n\to\infty }\EE \left(R^I_1(n)\right)\to0.
\end{align*}
Analogously, we also have as $n\to\infty$,
\[
\lim_{n\to\infty }\EE\left( R^O_1(n)\right)  \to 0.
\]
Hence, $R_1(n)\convp 0$ as $n\to\infty$.

To prove $R_2(n)\convas 0$, first consider for $k\ge 1$, 
\begin{align*}
W_k& \equiv B_k
\ind_{\left\{\left(BI^{(|\barbV(k)|)}_{\deltain}\left(\frac{1}{1+\deltain p}\log\left(\frac{n}{k}\right)\right), 1+ \widetilde{BI}^{(|\barbV(k)|)}_{1+\deltaout}\left(\frac{1-p}{1+\deltaout p}\log\left(\frac{n}{k}\right)\right)\right) = (m,l)\right\}} \\
& - p\PP\left(\left(BI^{(|\barbV(k)|)}_{\deltain}\left(\frac{1}{1+\deltain p}\log\left(\frac{n}{k}\right)\right), 1+ \widetilde{BI}^{(|\barbV(k)|)}_{1+\deltaout}\left(\frac{1-p}{1+\deltaout p}\log\left(\frac{n}{k}\right)\right)\right) = (m,l)\right).
\end{align*}
Note that $\EE(W_k) =0$ for all $k\ge 1$. Also, for $k> j\ge 1$, if $B_k=1$, then
$|\barbV(k)|>|\barbV(j)|$. 
Define in addition $${W}^*_k:= 
\ind_{\left\{\left(BI^{(|\barbV(k)|)}_{\deltain}\left(\frac{1}{1+\deltain p}\log\left(\frac{n}{k}\right)\right), 1+ \widetilde{BI}^{(|\barbV(k)|)}_{1+\deltaout}\left(\frac{1-p}{1+\deltaout p}\log\left(\frac{n}{k}\right)\right)\right) = (m,l)\right\}},$$
and we have
\begin{align*}
\PP&\left(B_k{W}^*_k=1, B_j{W}^*_j=1\right)
= \EE\left(\ind_{\{B_k=1\}}\ind_{\{B_j=1\}}
\PP^{\{B_i\}_{i=1}^k}\left({W}^*_k=1, {W}^*_j=1\right)
\right),\\
\intertext{where applying the independence among $\{BI_{\deltain}^{(v)}(\cdot)\}_{v\ge 2}$, $\{\widetilde{BI}_{1+\deltaout}^{(v)}(\cdot)\}_{v\ge 1}$, and $\{B_k:k\ge 1\}$, implies}
&= \PP\left(B_k{W}^*_k=1\right) \PP\left(B_j{W}^*_j=1\right).
\end{align*}
Hence, in general, for $k\neq j$, we have $\EE(B_kW^*_k B_jW^*_j) = \EE(B_kW^*_k) \EE(B_jW^*_j)$,
then
$\EE(W_k^3 W_j) = \EE(W_k^2 W_j) = \EE(W_k W_j) = 0$.
Then by the Markov's inequality, for any $\epsilon>0$,
\begin{align*}
\PP&\left(\left|\frac{1}{n}\sum_{k=1}^n W_k\right|\ge \epsilon\right) \le \frac{1}{n^4 \epsilon^4} \EE \left(\sum_{k=1}^nW_k\right)^4\\
&= \frac{1}{n^4 \epsilon^4} \EE\left(\sum_{k=1}^n W_k^4 + 4\sum_{k\neq l} W_k W_l^3
+ 3\sum_{k\neq l} W_k^2 W_l^2
+ 6 \sum_{k\neq l\neq i} W_k^2 W_l W_i
+ \sum_{k\neq l\neq i\neq j} W_k W_l W_i W_j
\right)\\
& = \frac{1}{n^4 \epsilon^4} \EE\left(\sum_{k=1}^n W_k^4 
+ 3\sum_{k\neq l} W_k^2 W_l^2
\right),\\
\intertext{since $|W_k|\le 1$ for $k\ge 1$, we have}
&\le \frac{1}{n^3\epsilon^4} + \frac{3}{n^2\epsilon^2}.
\end{align*}
Then by Borel-Cantelli lemma, we have
\beqq\label{eq:conv_Wk}
\frac{1}{n}\sum_{k=1}^n W_k \convas 0.
\eeqq
Since $|\barbV(n)|/n\convas p$, \eqref{eq:conv_Wk} implies
\[
R_2(n)\convas 0,\qquad\text{as }n\to\infty.
\]

By \cite[Equation~(2.2)]{tavare:1987},
 a birth immigration process $BI_\delta(\cdot)$ satisfies for fixed $t\ge 0$,
\[
BI_\delta(t) \stackrel{d}{=} Z_\delta(e^{-t}).
\]
Since the pmf of $Z_\delta(e^{-t})$ is bounded and continuous in $t$, then by the Riemann integrability of 
$$ \PP\left(\left(Z_{\deltain}\left(t^{1/(1+\deltain p)}\right), 1+ \widetilde{Z}_{1+\deltaout}\left(t^{(1-p)/(1+\deltaout p)}\right)\right) = (m,l)\right),$$
 we see $R_3(n)\convas 0$ as $n\to\infty$,
which completes the proof of the corollary. 
\end{proof}

\subsection{Proofs in Section~\ref{sec:POPA}}\label{subsec:pf2}
\begin{proof}[Proof of Proposition~\ref{prop:conv_deg_PGPA}]
Define the sigma-algebra
\[
\mathcal{G}_n := \sigma\left\{\bigl(\bfI(k), \bfO(k)\bigr): 1\le k\le n\right\}.
\]
By the definition of $\bfI(\cdot)$, we have for $v\in \bfV(n)$,
\begin{align*}
\EE^{\mathcal{G}_n}\left(I_v(n+1)\right) &= I_v(n) + \frac{(\lambda+1)\left(I_v(n)+\deltain\right)}{M_n+1+\deltain |\bfV(n)|}.
\end{align*}
Therefore, for $v\in\bfV(i)\setminus\bfV(i-1)$, $i\ge 0$,
\[
\mathcal{M}_n := \frac{I_v(n)+\deltain}{\prod_{k=0}^{n-1}\left(1+\frac{\lambda+1}{M_k+1+\deltain |\bfV(k)|}\right)}
\]
is a non-negative $\mathcal{G}_n$-martingale. By the martingale convergence theorem, $\mathcal{M}_n$ converges a.s. to
some limit $\xi_i$. 

Since for $x\ge 0$, we have $\frac{x}{1+x}\le \log(1+x)\le x$, then
\begin{align}
\exp\left\{\sum_{k=0}^{n-1}\frac{\lambda+1}{M_k+1+\deltain |\bfV(k)|+\lambda+1}\right\}&\le
\prod_{k=0}^{n-1} \left(1+\frac{\lambda+1}{M_k+1+\deltain |\bfV(k)|}\right)\nonumber\\
& \le \exp\left\{\sum_{k=0}^{n-1}\frac{\lambda+1}{M_k+1+\deltain |\bfV(k)|}\right\}
\label{eq:order}
\end{align}
By \cite[Theorem III.9.4, Page 120]{athreya:ney:2004},
we have that there exists some limiting random variable $Y'$ such that
\[
\sum_{k=1}^{n-1}\left( \frac{\lambda+1}{M_k+1+ \deltain |\bfV(k)|} - \frac{1}{(1+\deltain p)k} \right)
\convas Y'.
\]
Since $\sum_{k=1}^{n-1}1/k -\log n\to\gamma$, which is Euler's constant, then
\[
\sum_{k=1}^{n-1}\frac{\lambda+1}{M_k+ 1+\deltain |\bfV(k)|} - \frac{1}{1+\deltain p}\log(n)
\convas Y'+\frac{\gamma}{1+\deltain p} =: Y_U.
\]
Similarly, there exists another random variable $Y_L$ such that
\[
\sum_{k=1}^{n-1}\frac{\lambda+1}{M_k+1+ \deltain |\bfV(k)|+\lambda+1} - \frac{1}{1+\deltain p}\log(n)
\convas Y_L.
\]
Therefore, we have
\begin{align*}
\limsup_{n\to\infty}&\,\frac{1}{n^{1/(1+\deltain p)}}\prod_{k=0}^{n-1}\left(1+\frac{\lambda+1}{M_k+1+\deltain |\bfV(k)|}\right)\le e^{Y_U}<\infty,\\
\limsup_{n\to\infty}&\,{n^{1/(1+\deltain p)}}\prod_{k=0}^{n-1}\left(1+\frac{\lambda+1}{M_k+1+\deltain |\bfV(k)|}\right)^{-1}\le e^{-Y_L}<\infty.
\end{align*}


The proof machinery above is also applicable to $O_v(n)$, which then completes the proof of the proposition.
\end{proof}

\begin{proof}[Proof of Theorem~\ref{thm:IO_Pois}]
Note that for $v\in \bfV(i)\setminus\bfV(i-1)$, $n\epsilon \le i\le n$, $\epsilon>0$,
\begin{align*}
\PP^{\mathcal{G}_{i}}&\left(I_v(i+1) = 0\right)\\
&= \left(1-\frac{\deltain}{M_{i}+1+\deltain |\bfV(i)|}\right) e^{-\frac{\lambda \deltain}{M_{i}+1+\deltain |\bfV(i)|}}\\
&= \left(1-\frac{\deltain}{M_{i}+1+\deltain |\bfV(i)|}\right) \left(1-\frac{\lambda\deltain}{M_{i}+1+\deltain |\bfV(i)|} + O(M_i^{-2})\right)\\
&= 1-\frac{(\lambda+1)\deltain}{M_{i}+1+\deltain |\bfV(i)|} + O(M_i^{-2});\\
\intertext{since as $i\to\infty$,
$
\left(\frac{(\lambda+1)\deltain}{M_{i}+1+\deltain |\bfV(i)|}\right)/\left(\frac{\deltain}{i+1+\deltain |\barbV(i)|}\right)
\convas 1$, then we can compare the Poisson PA model with the traditional one:}
&= 1- \frac{\deltain}{i+1+\deltain |\barbV(i)|} (1+ o(1)) + O(M_i^{-2})\\
&= \PP^{\mathcal{F}_{2i}} \left(\barI_{|\barbV(i)|}(2(i+1)) = 0\right) + o(M_i^{-1}),
\end{align*}
\sid{by \eqref{eq:J_PA}}.
For $m\ge 1$,
\begin{align}
\PP^{\mathcal{G}_{i}}&\left(I_v(i+1) = m\right)\nonumber \\
=& \frac{\deltain}{M_{i}+1+\deltain |\bfV(i)|}\left(\frac{\lambda\deltain}{M_{i}+1+\deltain |\bfV(i)|}\right)^{m-1}\frac{e^{-\frac{\lambda \deltain}{M_{i}+1+\deltain |\bfV(i)|}}}{(m-1)!}\nonumber\\
&+ \left(1-\frac{\deltain}{M_{i}+1+\deltain |\bfV(i)|}\right)\left(\frac{\lambda\deltain}{M_{i}+1+\deltain |\bfV(i)|}\right)^{m}\frac{e^{-\frac{\lambda \deltain}{M_{i}+1+\deltain |\bfV(i)|}}}{m!}.
\label{eq:approx_in}
\end{align}
If $m=1$, then
\begin{align*}
\PP^{\mathcal{G}_{i}}&\left(I_v(i+1) = 1\right)\\
=&\left( \frac{\deltain}{M_{i}+1+\deltain |\bfV(i)|}+\frac{\lambda\deltain}{M_{i}+1+\deltain |\bfV(i)|}\left(1-\frac{\deltain}{M_{i}+1+\deltain |\bfV(i)|}\right)\right)\\
&\qquad\times \exp\left\{-\frac{\lambda \deltain}{M_{i}+1+\deltain |\bfV(i)|}\right\}\\
=& \left( \frac{(\lambda+1)\deltain}{M_{i}+1+\deltain |\bfV(i)|} + O(M_i^{-2})\right)\left(1-\frac{\lambda \deltain}{M_{i}+1+\deltain |\bfV(i)|} + O(M_i^{-2})\right)\\
=& \frac{(\lambda+1)\deltain}{M_{i}+1+\deltain |\bfV(i)|} + O(M_i^{-2})\\
=& \PP^{\mathcal{F}_{2i}} \left(\barI_{|\barbV(i)|}(2(i+1)) = 1\right) +o(M_i^{-1}).
\end{align*}
For $m\ge 2$, 
\eqref{eq:approx_in} implies 
$\PP^{\mathcal{G}_{i}}\left(I_v(i+1) = m\right) = O(M_i^{-2})$.
A similar argument also applies to $O_v(i+1)$, and we have
for $n\epsilon \le i\le n$,
\begin{align*}
\PP^{\mathcal{G}_{i}}&(I_v(i+1)=m, O_v(i+1)=l)\\
 &= \PP^{\mathcal{G}_{i}}(I_v(i+1)=m)\PP^{\mathcal{G}_{i}}( O_v(i+1)=l)\\
&= \left(\PP^{\mathcal{F}_{2i}}(\barI_{|\barbV(i)|}(2(i+1))=m)+o(M_i^{-1})\right) \left(\PP^{\mathcal{F}_{2i}}(\barO_{|\barbV(i)|}(2(i+1))=l) +o(M_i^{-1})\right)\\
&= \PP^{\mathcal{F}_{2i}}(\barI_{|\barbV(i)|}(2(i+1))=m, \barO_{|\barbV(i)|}(2(i+1))=l) +o(M_i^{-1}).
\end{align*}

Similarly, the conditional distribution $I_v(i+2)$ given $\mathcal{G}_{i+1}$ satisfies
\begin{align*}
\PP^{\mathcal{G}_{i+1}}&\left(I_v(i+2) = m-I_v(i+1) \right)\\
&=  \PP^{\mathcal{F}_{2(i+1)}}\left(\barI_{|\barbV(i)|}(2(i+2)) = m-\barI_{|\barbV(i)|}(2(i+1)) \right) + o(M_{i+1}^{-1}).
\end{align*}
Since 
\begin{align*}
\PP^{\mathcal{G}_{i}}&\left(I_v(i+2) = m\right)\\
&= \sum_{u=0}^m \PP^{\mathcal{G}_{i}}\left(I_v(i+2) = m\middle| I_v(i+1)=u, M_{i+1},
\{B_k\}_{k=M_i+1}^{M_{i+1}}\right)\PP^{\mathcal{G}_{i}}\left(I_v(i+1) = u\right),
\end{align*}
we then have
\begin{align*}
\PP^{\mathcal{G}_{i}}\left(I_v(i+2) = m\right)
&= \PP^{\mathcal{F}_{2i}}\left(\barI_{|\barbV(i)|}(2(i+2)) = m \right) + o(M_i^{-1}),
\end{align*}
which leads to 
\begin{align*}
\PP^{\mathcal{G}_{i}}&\left(I_v(i+2) = m, O_v(i+2)=l\right)\\
&= \PP^{\mathcal{F}_{2i}}\left(\barI_{|\barbV(i)|}(2(i+2)) = m,\barO_{|\barbV(i)|}(2(i+2))=l \right) + o(M_i^{-1}).
\end{align*}
Using such induction steps to proceed, we conclude that for $n\epsilon\le i\le n$, $\epsilon>0$,
\begin{align}
\PP^{\mathcal{G}_{i}}\left(I_v(n) = m, O_v(n)=l\right)
&= \PP^{\mathcal{F}_{2i}}\left(\barI_{|\barbV(i)|}(2n) = m,\barO_{|\barbV(i)|}(2n)=l \right) + o(M_i^{-1}).
\label{eq:approx_cond}
\end{align}
From \cite[Theorem III.9.4, Page 120]{athreya:ney:2004}, we have
$\EE(i/M_i)\to 1/(\lambda+1)$ as $i\to\infty$,
so taking the expectation on both sides of \eqref{eq:approx_cond} gives
\begin{align*}
\PP\left(I_v(n) = m, O_v(n)=l\right)
&= \PP\left(\barI_{|\barbV(i)|}(2n) = m,\barO_{|\barbV(i)|}(2n)=l \right) + o(i^{-1}),
\end{align*}
and after summing over $n\epsilon\le i\le n$, we have: there exists some constant $C>0$ (depending on $\deltain,\deltaout, p$) such that
\begin{align*}
\left|\frac{1}{n} \right. &\left.\sum_{i=n\epsilon}^n \PP\left(I_{|\bfV(i)|}(n) = m, O_{|\bfV(i)|}(n)=l\right)\right.\\
&\left.  -\frac{1}{n}\sum_{i=n\epsilon}^n \PP\left(\barI_{|\barbV(i)|}(2n) = m,\barO_{|\barbV(i)|}(2n)=l \right)\right|\le \frac{C\log n}{n}\stackrel{n\to\infty}{\longrightarrow} 0.
\end{align*}
Then \eqref{eq:approx_pij} follows.

To prove \eqref{eq:POPA_pij}, we first note that following a similar reasoning as in Corollary~\ref{cor:PA_count}, we can embed $\bfI(\cdot)$ into a sequence of iid special birth-immigration processes $\{\mathcal{B}^{(v)}_{\ind_{\{v=1\}}+\deltain}(t): t\ge0\}_{v\ge 1}$ as follows. 
After observing $M_1$ and $B_1,\ldots, B_{M_1}$, let $\mathcal{T}_1, \ldots, \mathcal{T}_{M_1}$, be the first $M_1$ jump times
of the $\mathcal{B}^{(1)}_{1+\deltain}(\cdot)$ process such that $\mathcal{B}^{(1)}_{1+\deltain}(\cdot)$ has a constant transition rate, $1+\deltain$, over $[0,\mathcal{T}_{M_1}]$.
For $k=1,\ldots, M_1$, if $B_k = 1$, then we start a new special birth-immigration process
$\left\{\mathcal{B}^{(1+\sum_{i=1}^k B_i)}_{\deltain}(t-\mathcal{T}_{M_1}): t\ge \mathcal{T}_{M_1}\right\}$.


For $n\ge 1$, given $\{M_k: k=1,\ldots, n+1\}$,
$\{B_k: k=1,\ldots, M_{n+1}\}$,
 and
\beqq\label{eq:BI_POPAin}
\left\{\mathcal{B}^{(v)}_{\ind_{\{v=1\}}+\deltain}(t-\mathcal{T}_{M_i}): t\ge \mathcal{T}_{M_i}\right\},
\quad v\in \bfV(i)\setminus\bfV(i-1),0\le i\le n,
\eeqq
we set $\mathcal{T}_{M_n+1}, \ldots, \mathcal{T}_{M_{n+1}}$ to be the next $\Delta M_{n+1}$ jump times among the $|\bfV(n)|$ processes in \eqref{eq:BI_POPAin},
such that $\mathcal{B}^{(v)}_{\ind_{\{v=1\}}+\deltain}(\cdot-\mathcal{T}_{M_i})$ has a
constant transition rate, $\mathcal{B}^{(v)}_{\ind_{\{v=1\}}+\deltain}(\mathcal{T}_{M_n}-\mathcal{T}_{M_i})+\ind_{\{v=1\}}+\deltain$, over $(\mathcal{T}_{M_n},\mathcal{T}_{M_{n+1}}]$.
For $k=M_n+1,\ldots, M_{n+1}$, use $\mathcal{J}_{k}'$, to denote which process in \eqref{eq:BI_POPAin} jumps at $\mathcal{T}_k$,
and if $B_k = 1$, start a new special birth-immigration process 
$\left\{\mathcal{B}^{(1+\sum_{i=1}^k B_{i})}_{\deltain}(t-\mathcal{T}_{M_{n+1}}): t\ge \mathcal{T}_{M_{n+1}}\right\}$.
Then similar to the embedding results in Theorem~\ref{thm:PA_BI}, 
assuming $\mathcal{T}_0:=0$,
we have on $\mathbb{N}^\infty$,
\begin{align*}
\left\{\bfI({n}): n\ge 0\right\}
\stackrel{d}{=}& \left\{\left(\left\{\ind_{\{v=1\}} + \mathcal{B}^{(v)}_{\ind_{\{v=1\}}+\deltain}\left(\mathcal{T}_{M_n} - \mathcal{T}_{M_i}\right)\right\}_{v\in \bfV(i)\setminus\bfV(i-1),0\le i\le n}, 0,\ldots\right): n\ge 0\right\}.
\end{align*}

Following the embedding analogy above and as in Section~\ref{subsubsec:embed}, we also embed
the out-degree sequence,
$\left\{\bfO({n}): n\ge 0\right\}$, into a sequence of iid special birth-immigration processes 
$\{\widetilde{\mathcal{B}}^{(v)}_{1+\deltaout}(t): t\ge0\}_{v\ge 1}$.
With $M_1$ and $B_1,\ldots, B_{M_1}$ given, let $\widetilde{\mathcal{T}}_0:=0$ and $\widetilde{\mathcal{T}}_1, \ldots, \widetilde{\mathcal{T}}_{M_1+1-|\bfV(1)|}$, with $|\bfV(1)|=1+\sum_{k=1}^{M_1}B_k$, be the first $M_1+1-|\bfV(1)| = \sum_{k=1}^{M_1}(1-B_k)$ jump times
of the $\widetilde{\mathcal{B}}^{(1)}_{1+\deltaout}(\cdot)$ process such that
 $\widetilde{\mathcal{B}}^{(1)}_{1+\deltaout}(\cdot)$ has a constant transition rate, $1+\deltaout$, over $[0,\widetilde{\mathcal{T}}_{M_1+1-|\bfV(1)|}]$.
For $k=1,\ldots, M_1$, if $B_k = 1$, then we start a new special birth-immigration process
$$\left\{\widetilde{\mathcal{B}}^{(1+\sum_{i=1}^k B_i)}_{1+\deltaout}
\left(t-\widetilde{\mathcal{T}}_{M_1+1-|\bfV(1)|}\right): t\ge \widetilde{\mathcal{T}}_{M_1+1-|\bfV(1)|}\right\},$$
at time $\widetilde{\mathcal{T}}_{M_1+1-|\bfV(1)|}$.

For $n\ge 1$, conditional on $\{M_k: k=1,\ldots, n+1\}$,
$\{B_k: k=1,\ldots, M_{n+1}\}$, 
$|\bfV(i)| = 1+\sum_{k=1}^{M_i}B_k$, $0\le i\le n$,
and 
\beqq\label{eq:BI_POPAout}
\left\{\widetilde{\mathcal{B}}^{(v)}_{1+\deltaout}\left(t-\widetilde{\mathcal{T}}_{M_i+1-|\bfV(i)|}\right): 
t\ge \widetilde{\mathcal{T}}_{M_i+1-|\bfV(i)|}\right\},\quad
v\in \bfV(i)\setminus\bfV(i-1),0\le i\le n,
\eeqq
set $\widetilde{\mathcal{T}}_{M_n+2-|\bfV(n)|}, \ldots, \widetilde{\mathcal{T}}_{M_{n+1}+1-|\bfV(n+1)|}$ to be the next $\sum_{k=1}^{\Delta M_{n+1}} (1-B_{M_n+k})$ jump times among the $|\bfV(n)|$ processes in \eqref{eq:BI_POPAout},
such that $\widetilde{\mathcal{B}}^{(v)}_{1+\deltaout}(\cdot-\widetilde{\mathcal{T}}_{M_i+1-|\bfV(i)|})$ has 
a constant transition rate, $\widetilde{\mathcal{B}}^{(v)}_{1+\deltaout}(\widetilde{\mathcal{T}}_{M_n+1-|\bfV(n)|}-\widetilde{\mathcal{T}}_{M_i+1-|\bfV(i)|})+1+\deltaout$, over $\left(\widetilde{\mathcal{T}}_{M_n+1-|\bfV(n)|},\widetilde{\mathcal{T}}_{M_{n+1}+1-|\bfV(n+1)|}\right]$.
For $k=M_n+1,\ldots, M_{n+1}$, if $B_k=0$, use $\widetilde{\mathcal{J}}_{k}'$ to denote which process in \eqref{eq:BI_POPAout} jumps at time $\widetilde{\mathcal{T}}_{\sum_{i=1}^k (1-B_{i})}$.
If $B_k = 1$, we set $\widetilde{\mathcal{J}}_{k}'=0$ and start a new special birth-immigration process 
$$\left\{\widetilde{\mathcal{B}}^{(1+\sum_{i=1}^k B_{i})}_{1+\deltaout}(t-\widetilde{\mathcal{T}}_{M_{n+1}+1-|\bfV(n+1)|}): t\ge \widetilde{\mathcal{T}}_{M_{n+1}+1-|\bfV(n+1)|}\right\}.$$ 
Then we have on $\left(\mathbb{N}_{>0}\right)^\infty$,
\begin{align*}
&\left\{\bfO({n}): n\ge 0\right\}\\
\stackrel{d}{=}& \left\{\left(\left\{1 + \widetilde{\mathcal{B}}^{(v)}_{1+\deltaout}\left(\widetilde{\mathcal{T}}_{M_n+1-|\bfV(n)|} - \widetilde{\mathcal{T}}_{M_i+1-|\bfV(i)|}\right)\right\}_{v\in \bfV(i)\setminus\bfV(i-1),0\le i\le n}, 0,\ldots\right): n\ge 0\right\}.
\end{align*}
By the embedding framework, we see that
\[
\mathcal{T}_{M_n} - \sum_{k=0}^{n-1} \frac{\lambda+1}{M_k+1+\deltain |\bfV(k)|},\qquad\text{and}\qquad \widetilde{\mathcal{T}}_{M_n+1-|\bfV(n)|} - \sum_{k=0}^{n-1} \frac{(\lambda+1)(1-p)}{M_k+1+\deltaout |\bfV(k)|}
\]
are $L_2$-bounded martingales with respect to sigma fields
\begin{align*}
&\sigma\left(\left\{M_k\right\}_{k=1}^n;\left\{B_k\right\}_{k=1}^{M_n} ;\left\{{\mathcal{B}}^{(v)}_{\ind_{\{v=1\}}+\deltain}(t-{\mathcal{T}}_{M_i}): {\mathcal{T}}_{M_i}\le t\le {\mathcal{T}}_{M_n}\right\}_{v\in\bfV(i)\setminus\bfV(i-1)}, 0\le i\le n\right),\\
&\sigma\left(\left\{M_k\right\}_{k=1}^n;\left\{B_k\right\}_{k=1}^{M_n} ;\right.\\
&\left.\quad
\left\{\widetilde{\mathcal{B}}^{(v)}_{1+\deltaout}(t-\widetilde{\mathcal{T}}_{M_i+1-|\bfV(i)|}): \widetilde{\mathcal{T}}_{M_i+1-|\bfV(i)|}\le t\le \widetilde{\mathcal{T}}_{M_n+1-|\bfV(n)|}\right\}_{v\in\bfV(i)\setminus\bfV(i-1)}, 0\le i\le n\right),
\end{align*}
respectively.

Let $\{\mathcal{B}_{\deltain}(t): t\ge 0\}$ and $\{\widetilde{\mathcal{B}}_{1+\deltaout}(t): t\ge 0\}$
be two independent special birth-immigration processes
which have the same behavior as $\{\mathcal{B}^{(2)}_{\deltain}(t):t\ge0\}$ and
$\{\widetilde{\mathcal{B}}^{(2)}_{1+\deltaout}(t):t\ge 0\}$, respectively.
Now we follow the reasoning in the proof of Corollary~\ref{cor:PA_count} to obtain
\begin{align*}
\frac{1}{|\bfV(n)|}&\sum_{v=1}^{|\bfV(n)|} \ind_{\{I_{v}(n)=m, O_{v}(n)=l\}}\\
&\convp \int_0^1 \PP\left(\mathcal{B}_{\deltain}\left(-\frac{1}{1+\deltain p}\log t\right)=m,
1+\widetilde{\mathcal{B}}_{1+\deltaout}\left(-\frac{1-p}{1+\deltaout p}\log t\right)=l\right)\dd t\\
&=: p'_{m,l}.
\end{align*}
Since a.s.
\[
\frac{1}{|\bfV(n)|}\sum_{v=1}^{|\bfV(n)|} \ind_{\{I_{v}(n)=m, O_{v}(n)=l\}}\le 1,
\]
then as $n\to\infty$,
\begin{align*}
\EE\left(\frac{1}{|\bfV(n)|}\sum_{v=1}^{|\bfV(n)|} \ind_{\{I_{v}(n)=m, O_{v}(n)=l\}}\right) 
&\longrightarrow p'_{m,l}.
\end{align*}
By the approximation in \eqref{eq:approx_pij}, we see that
\[
p'_{m,l} = p_{m,l},
\]
which completes the proof of \eqref{eq:POPA_pij}.
\end{proof}

\bibliographystyle{imsart-number} 
\bibliography{./bibfile_0512.bib}
\end{document}